\documentclass[11pt]{article}

\usepackage[normalem]{ulem} 
\usepackage{cancel}

\usepackage{graphicx,amsmath,amsfonts,amssymb,color,bbm,bm}

\graphicspath{{figures/}}
\usepackage{subcaption}

\usepackage{hyperref}

\usepackage{tikz}

\usepackage{epsfig}

 \newcommand{\pend}{\hfill \thicklines \framebox(5.5,5.5)[l]{}}
 \newenvironment{proof}{\noindent {\sc  Proof.} \rm}{\pend}

\numberwithin{equation}{section}

 \newtheorem{theorem}{Theorem}

 \newtheorem{remark}{Remark}[section]

 \newtheorem{corollary}{Corollary}[section]
 \newtheorem{definition}{Definition}[section]
 
 \newtheorem{property}{Properties}

\makeatletter
 \newif\if@borderstar
 \def\bordermatrix{\@ifnextchar*{%
 \@borderstartrue\@bordermatrix@i}{\@borderstarfalse\@bordermatrix@i*}%
 }
 \def\@bordermatrix@i*{\@ifnextchar[{\@bordermatrix@ii}{\@bordermatrix@ii[()]}}
 \def\@bordermatrix@ii[#1]#2{%
 \begingroup
 \m@th\@tempdima8.75\p@\setbox\z@\vbox{%
 \def\cr{\crcr\noalign{\kern 2\p@\global\let\cr\endline }}%
 \ialign {$##$\hfil\kern 2\p@\kern\@tempdima & \thinspace %
 \hfil $##$\hfil && \quad\hfil $##$\hfil\crcr\omit\strut %
 \hfil\crcr\noalign{\kern -\baselineskip}#2\crcr\omit %
 \strut\cr}}%
 \setbox\tw@\vbox{\unvcopy\z@\global\setbox\@ne\lastbox}%
 \setbox\tw@\hbox{\unhbox\@ne\unskip\global\setbox\@ne\lastbox}%
 \setbox\tw@\hbox{%
 $\kern\wd\@ne\kern -\@tempdima\left\@firstoftwo#1%
 \if@borderstar\kern2pt\else\kern -\wd\@ne\fi%
 \global\setbox\@ne\vbox{\box\@ne\if@borderstar\else\kern 2\p@\fi}%
 \vcenter{\if@borderstar\else\kern -\ht\@ne\fi%
 \unvbox\z@\kern-\if@borderstar2\fi\baselineskip}%
 \if@borderstar\kern-2\@tempdima\kern2\p@\else\,\fi\right\@secondoftwo#1 $%
 }\null \;\vbox{\kern\ht\@ne\box\tw@}%
 \endgroup
 }
 \makeatother

\begin{document}

\title{Markov Modelling Approach for Queues with Correlated Service Times --- the $M/M_D/2$ Model}
\author{Suman Thapa\footnote{School of Mathematics and Statistics, Carleton University, Ottawa ON Canada K1S 5B6, email address: \newline \href{mailto:sumanthapa@cmail.carleton.ca}{sumanthapa@cmail.carleton.ca}} and Yiqiang Q. Zhao\footnote{School of Mathematics and Statistics, Carleton University, Ottawa ON Canada K1S 5B6, email address: \href{mailto:zhao@math.carleton.ca}{zhao@math.carleton.ca}} 
}
\date{(Updated, April 2026)}

\maketitle

\begin{abstract}
  Demand for studying queueing systems with multiple servers providing correlated services was created about 60 years ago, motivated by various applications. In recent years, the importance of such studies has been significantly increased, supported by new applications of greater significance to much larger scaled industry, and the whole society. Such studies have been considered very challenging.
  In this paper, a new Markov modelling approach for queueing systems with servers providing correlated services is proposed. We apply this new proposed approach to a queueing system with arrivals according to a Poisson process and two positive correlated exponential servers, referred to as the $M/M_D/2$ queue. We first prove that the queueing process (the number of customers in the system) is a Markov chain, and then provide an analytic solution for the stationary distribution of the process, based on which it becomes much easier to see the impact of the dependence on system performance compared to the performance with independent services. 
\end{abstract}

\section{Introduction}

\subsection{Motivations}

The $M/M/c$ queueing system was first introduced and studied by Agner Krarup Erlang in 1909 in \cite{Erlang:1909}, which later becomes one of the most popular queueing systems in applications with independent servers. Since the late 1970s, queueing systems with correlated servers have been proposed and studied by several researchers, motivated by various applications. 

Kleinrock realized the importance of dependent service times in message-switching communication networks (MSCN) in \cite{Kleinrock:1964}. In MSCN, messages (customers) are divided into pieces according to channel capacity (service time) during message transmission through the network, which requires identical service times at each tandem node in queueing language. This type of network models (or tandem queues) later becomes more popular under the packet-switching technique, and is the main focus of the studies for queueing systems with dependent servers.  

Mitchell \textit{et al.} in \cite{Mitchell et al:1977} offered two scenarios in production lines and logistics, in which service times could be correlated: (1) in a paper mill, and (2) in the overhaul of an aircraft engine. They suggested that ``An obvious approach is to use a multivariate exponential distribution with non-zero correlations in place of the usual independent exponential service times.'' And they further questioned that ``it is not clear that the birthdeath equation approach can be modified to incorporate dependent service times. Moreover, any such formulation would very likely be analytically intractable.''

In \cite{Choo-Conolly:1980}, Choo and Conolly described two situations, in which correlated service times are intuitively obvious: ``That correlation may be relevant is intuitively clear from the supermarket example
where most customers who spend a long time in the shopping area also ipso facto require a long checkout time.
And a patient with an unknown disease may spend a long time at each of a series of investigative stages until the
disease is identified, followed perhaps by long periods of treatment and convalescence.'' 

Two more applications of queues requiring dependent service times were provided in Pang and Whitt \cite{Pang-Whitt:2012-MSOM}. The first example is in a technical support telephone call center responding to service call. It is anticipated that the handling (service) times of calls, after the product with a defect is first introduced, are longer-than-usual. The other example is in a hospital emergency room.
Multiple patients may be associated with the same medical incident (say a highway accident, or food poisoning at the same restaurant), all of whom require longer service times.

In recent years, new applications, in which service times provided by different servers are correlated, lead to additional motivations of studies of queueing system with correlated servers, including this one. Our model, referred to as the $M/M_D/2$ queue, is a generalization of the standard $M/M/2$ queueing system by allowing correlations between the two servers. Specifically, the service times provided by the two servers follow the Marshall-Olkin bivariate exponential distribution (MO-BVED, and MO-BVE for Marshall-Olkin bivariate exponential). MO-BVED describe a type of positive correlation between the two service times and its unique feature is to allow simultaneous departures from both servers (see next section for more details). This model can be used in many applications, such as:
\begin{itemize}
  \item \textbf{Two servers in the same rack with shared power and cooling.}
    Simultaneous departures can be caused by rack--level power dip, top-of-rack switch failure, cooling failure, operating-system upgrades that reboot the whole machine.

  \item \textbf{Two parallel database replicas or storage devices.}
  In this case, two mirrored disks or database replicas handle requests in parallel with simultaneous departures due to controller failure, shared filesystem corruption, SAN outage, data-centre network partition.

  \item \textbf{Two human servers subject to common interruptions.} 
For example, two doctors in the same ward, or two agents in a small call centre, serve different customers. In this case, fire alarms, ward-level emergencies, IT outages that could force everyone to stop service.

  \item \textbf{Two communication links with a shared critical component.}
For example, two VPN tunnels or logical links traverse the same physical fibre or core router. In this case, fibre cuts, core-router reboots, or line-card failures can lead to simultaneous departures.

    \item \textbf{Queueing systems with negative customers (signals).} It is very interesting to note that the $M/M_D/2$ model and queueing systems with negative customers (signals), often referred to as $G$-networks proposed by Gelenbe~\cite{Gelenbe:1989,Gelenbe:1991}, share a same feature: simultaneous departures. Our model provides a tool for studying a new type of $G$-networks, where simultaneous departure rate is independent of single departure rates (in $G$-networks defined by Gelenbe, simultaneous departure rate and single departure rate are correlated). 
\end{itemize}

The $M/M_D/2$ model also has potential uses in applications, where service time distribution is not MO-BVE, but can be approximated by the MO-BVED, including an extension of the standard $M/M/2$ queue to the model with peer-review pressure. For example, consider two cashiers next to each other serving customers. Each cashier checks the status of the other one after a random time to make a decision of speedup or slowdown of its service due to peer-pressure. When the speedup/slowdown mechanism would cause simultaneous departures (or service completions at roughly the same time), the $M/M_D/2$ model could be a good approximation for the application.

\subsection{Literature review}

Though it is important and well-motivated, the number of literature studies on queueing systems with correlated service times is still relatively small.
It was believed, say by Choo and Conolly \cite{Choo-Conolly:1980} and Pinedo and Wolff~\cite{Pinedo-Wolff:1982}, that the work by Mitchell \textit{et al.} \cite{Mitchell et al:1977} in 1977 was the first publication dealing with queueing systems with dependent service times provided by two servers in a tandem queue setting. Since then, most of studies have focused on tandem queues due to the nature of their structures. 
We should note that in fact, Kelly \cite{Kelly:1976} in 1976 (before \cite{Mitchell et al:1977} in 1977) already obtained some exact results for a time-sharing queueing network in which the service times for a customer remain the same throughout the network.

In \cite{Mitchell et al:1977}, the authors provided a simulation study on the impact of system performance of a two-node tandem queue with Poisson arrivals and bivariate exponential service times. They showed that the system behaviour is quite sensitive to the dependence, especially at higher utilizations, for either positive or negative correlation. 
The same tandem queueing system was also studied by Hoon Choo and Conolly~\cite{Choo-Conolly:1980} for the positive correlated case and obtained the Laplace transform of the (stationary) waiting time at stage 2.
 

In Boxma~\cite{Boxma:1979,Boxma:1979b}, the author also studied a two-stage tandem queue with Poisson arrivals, in which the service time at node 1 was extended from an exponential random variable (r.v.) to a r.v. with a general distribution. 
Each customer's service time at node two is identical to that at node one.  
Explicit (not in transformation form) expressions for the stationary sojourn time and waiting time were obtained in \cite{Boxma:1979}.
Computational implementations were addressed in \cite{Boxma:1979b}. 
For the same model, Boxma and Deng~\cite{Boxma-Deng:2000} obtained asymptotic results for delay distributions in the case of regularly varying service times, and also in heavy traffic conditions. In the book \cite{Kelly:1979}, Kelly also considered networks of symmetric queues allowing dependent service times.  

In two IBM technical reports \cite{Calo:1979,Calo:1980}, Calo obtained an expression for the Laplace-Stieltjes transform for the stationary distribution of the total waiting time of a customer in a tandem queue of multiple stages with a general arrival and general, but equal, service times. 

Pinedo and Wolff~\cite{Pinedo-Wolff:1982} considered Markovian tandem queues of $m$ stages (Poisson arrivals and exponential service times), in which each given customer has the equal (random) service time at each stage. 
Properties on stochastic ordering were used for their studies.
The work on the throughput in \cite{Pinedo-Wolff:1982} for the blocking case was extended for the non-blocking case by Browning~\cite{Browning:1998}.

Kelly~\cite{Kelly:1982} considered a series of nodes (queues) with a finite buffer between two consecutive nodes for transmitting messages. The lengths of successive messages are assumed to be i.i.d. random variables. The author obtained the rate at which buffer sizes need to increase in order to maintain the system throughput as the number $n$ of nodes increases.

Light traffic asymptotics for the expected delay in a series of $r \geq 2$ queues with correlated service times was investigated in Wolff~\cite{Wolff:1982}. For $r=2, 3$, the study applies to an arbitrary correlated joint distribution of service times, but for $r>3$, conclusions were made for equal service times. 

Ziedins~\cite{Ziedins:1993} also considered a tandem queueing system with a finite waiting capacity between two consecutive nodes, where any given customer has the same service time at each of the nodes. The author showed that, in terms of service time distributions with only two support points, it is always optimal to allocate the capacity as uniformly as possible, even when blocking occurs, which is different from the earlier suggestions in the literature.

Avi-Itzhak and Levy \cite{Avi-Itzhak-Levy:2001} considered tandem queues with deterministically correlated service times.

For tandem queues and through simulations results,
in Sandmann \cite{Sandmann:2007}, different types of correlations, including equal service times, were considered for the two node case with exponential and uniform service times at the first node.
Later, Sandmann~\cite{Sandmann:2012} considered the end-to-end delay with correlated service times, which is a continuation of his earlier work~\cite{Sandmann:2010}. Specifically, the service time at the first node is a random variable, and the service time at any other node is correlated with that at the first node. 

Pang and Whitt explored heavy traffic approximations under a different setting of queue models with dependent service times. Specifically, (a)
in \cite{Pang-Whitt:2013} (and also in \cite{Pang-Whitt:2012-MSOM}), an infinite-server queue, denoted by $G_t/G^D/\infty$, was considered, where the arrival process was assumed to satisfy a functional central limit theorem (FCLT), and the successive service times were assumed to be weakly dependent in the sense that the dependence among
the service times is limited so that the CLT remains valid, but the variability constant in the CLT is affected by the cumulative correlations. For such a system, heavy-traffic limits, including a functional weak law of large
numbers (FWLLN) and an FCLT were obtained; and
(b) in \cite{Pang-Whitt:2012}, the work in \cite{Pang-Whitt:2013} was extended to a generalized model allowing batch arrivals, denoted by $G_t^B/G^D/\infty$, in which dependent service times are restricted for customers in the same batch.
In \cite{Pang-Whitt:2012}, two special dependent service times, multivariate Marshall–Olkin (MO) exponential distributions and multivariate MO hyperexponential distributions within a batch were adopted for illustrating general results.

\subsection{Main contributions in this paper}
    The main contribution made in this paper is the proposal of a new approach for studying queueing systems with     
correlated servers, say the $M/M_D/2$ model. This new approach addresses a long-standing open concern raised in \cite{Mitchell et al:1977} (see also Section 1.1 on motivations). With our new approach, the birth-death equation approach indeed can be extended to deal with queueing systems with dependent service times. To prove the queueing process is actually a Markov chain, special cares have to be taken, including the dependent service rate, and the simultaneous departures from both servers due to singularities in the service time (see Theorem~\ref{the:1} and Corollary~\ref{cor:1}). It is worthwhile to point out that the Markov chain is no longer a birth-and-death process as in the case of independent servers. However, the stationary distribution is still geometric (see Theorem~\ref{the:2} and Corollary~\ref{cor:2}), which ensures that our formulation is indeed analytically tractable.

\subsection{Organization of the paper}

Following the introduction section, the rest of the paper is organized as follows: 
    in Section~\ref{sec:2}, we provide some literature results, which are required in our study, including the 
concept of 2-dimensional lack-of-memory property, and the MO-BVED and its properties;
    in Section~\ref{sec:3}, we introduce the $M/M_D/2$ queueing system and define the queueing process of this 
system. We prove that the queueing process of the $M/M_D/2$ system is a Markov chain;
    in Section~\ref{sec:4}, we provide the stationary solution of the Markov chain for the $M/M_D/2$ system, which    
is geometric with modified solutions at the boundary states;
    in Section~\ref{sec:5}, based on the main results obtained in previous sections, we demonstrate the impact of     
the dependence between the servers on the system performance; 
    and in the last section, Section~\ref{sec:6}, concluding remarks are made.

\section{Bivariate dependent service times} \label{sec:2}

Recall that the focus of this paper is to propose a Markovian model to study multi-server queues with dependent service times, in a similar fashion as $M/M/c$ was proposed to study multi-server queues with independent service times. 
In the case of the standard $M/M/c$ queue, the number of customers in the system at time $t$ is a Markov chain, proved based on the lack-of-memory property of the exponential random variable (r.v.). A univariate non-negative continuous r.v. $X$ is called lack-of-memory, if it satisfies the following property:
\[
    P(X > t +x | X > x) = P(X > t), \quad \text{for all $x, t \geq 0$.}
\] 
It is well-known that the exponential r.v. is the unique continuous non-negative r.v. satisfying the above lack-of-memory property.

In the literature, extending the lack-of-memory concept to multi-variate r.v.s was started from the seminal work by Marshall and Olkin~\cite{Marshall-Olkin:1967a}, and has been since a central topic in multi-variate distributions. References on this topic are vast, and besides \cite{Marshall-Olkin:1967a} we only mention Kots \textit{et al.} \cite{Kotz-et at:2000} and Lin \textit{et al.} \cite{Lin-Dou-Kuriki:2019} since all we need in our analysis can be conveniently found in these references. The lack-of-memory concept, in both strong and weak senses, and their properties introduced in this section are all literature results.

For the case of two r.v.s, it is intuitive to extend the lack-of-memory property to:
\[
    P(X > x+s, Y > y+t | X >s, Y>t) =P(X>x, Y>y),
\]
which is equivalent to, in terms of the survival function $\bar{H}$,
\[
    \bar{H}(x+s,y+t) = \bar{H}(x,y) \bar{H}(s,t), \quad \text{all $x, y, s, t \geq 0$,}
\] 
where $H$ is the distribution function of $(X,Y)$ and 
\[
    \bar{H}(x,y) = 1 - H(x,y) - P( X \leq x, Y >y) - P(X>x, Y \leq y).
\]
However, the above equation has only one solution, in which $X$ and $Y$ are independent exponential r.v.s. 
Therefore, lack-of-memory in weak sense has been proposed: $(X,Y)$ satisfies the bivariate lack-of-memory (BLM) property if
\[
    \bar{H}(x+t,y+t) = \bar{H}(x,y) \bar{H}(t,t), \quad \text{all $x, y, t \geq 0$.}
\] 
Lack-of-memory in weak sense is equivalent to:
\[
    P(X > x+t, Y > y+t | X >t, Y>t) =P(X>x, Y>y),\quad \text{all $x, y, t \geq 0$,}
\]
(see, for example, equation (2.7) in \cite{Marshall-Olkin:1967a}), which is more often used for the context of reliability in the literature. Unfortunately, the above expression, in terms of a conditional probability, is not a convenient form in our analysis for proving that the queueing process is a Markov chain. Instead, we use another equivalent form:
\begin{equation}\label{BLM-weak}
    P(X > x+t, Y > y+t | X >x, Y>y) =P(X>t, Y>t),\quad \text{all $x, y, t \geq 0$,}
\end{equation}
(see, for example, equation (2.9) in \cite{Marshall-Olkin:1967a}).

In our case, to have a Markov chain model, it is required that the marginal distributions are exponential (dealing with the case when one server is idle). Under this requirement (exponential marginals), it is well-known that the Marshall-Olkin bivariate exponential distribution is the only solution, which satisfies the lack-of-memory property in weak sense. 
\begin{definition} A bivariate distribution is called MO-BVE with parameters $\mu_1 \geq 0$, $\mu_2\geq 0$ and $\mu_{12}\geq 0$, if its survival function $\bar{H}$ is given by
\begin{equation}\label{eqn:M-O}
  \bar{H}(x,y) = e^{-(\mu_1 x + \mu_2 y + \mu_{12} \max(x,y))}, \quad x, y \geq 0,
\end{equation}
\end{definition}
For the MO-BVED, we provide a summary of its properties, which are needed in our analysis. 
\begin{property}
For the bivariate r.v.s $(X,Y)$ having the MO-BVED, the following properties hold:
\begin{description}
  \item[(1)] The MO-BVED is the only distribution satisfying the BLM property in weak sense (see, for example, equation \eqref{BLM-weak}), if the two marginal distributions are required to be exponential;

  \item[(2)] The density function $f$ of MO-BVE r.v.s $(X,Y)$ (or variables $(X,Y)$ having the MO-BVED) is given by
\begin{equation}\label{eqn:density of MO-BED}
  f(x,y)= \left \{ \begin{array}{ll}
            \mu_2 (\mu_1+\mu_{12}) \bar{H}(x,y), & \mbox{for } 0<y<x, \\
            \mu_1 (\mu_2+\mu_{12}) \bar{H}(x,y), & \mbox{for } 0<x<y, \\
            \mu_{12} \bar{H}(z,z), & \mbox{for } x=y=z >0,
                   \end{array} \right.
\end{equation}
or
\begin{equation}\label{eqn:density of MO-BED-b}
  f(x,y)= \left \{ \begin{array}{ll}
            \mu_2 (\mu_1+\mu_{12}) e^{-(\mu_1+\mu_{12})x +\mu_2 y}, & \mbox{for } 0<y<x, \\
            \mu_1 (\mu_2+\mu_{12}) e^{-(\mu_2+\mu_{12})y +\mu_1 x}, & \mbox{for } 0<x<y, \\
            \mu_{12} e^{-(\mu_1+\mu_2+\mu_{12}) z}, & \mbox{for } x=y=z >0.
                   \end{array} \right.
\end{equation}
See, for example, equation (47.51) in \cite{Kotz-et at:2000}. It is worthwhile to note that the MO-BVED is singular on $X=Y$ since $P(X=Y) >0$;

    \item[(3)] The marginal distributions of $X$ and $Y$ are both exponential with rates $\mu_1 + \mu_{12}$ and $\mu_2 + \mu_{12}$, respectively:
\begin{equation}\label{eqn:marginals}
   X \sim \mathrm{Exp}(\mu_1+\mu_{12}); \quad Y \sim \mathrm{Exp}(\mu_2+\mu_{12});
\end{equation}
    
    \item[(4)] Based on (\ref{eqn:density of MO-BED-b}), we have
\begin{equation}\label{eqn:X=Y>t}
  P(X=Y>t) = \frac{\mu_{12}}{\mu_1+\mu_2+\mu_{12}} e^{-(\mu_1+\mu_2+\mu_{12})t}, \quad t \geq 0,
\end{equation}
and
\begin{equation}\label{eqn:X=Y}
  P(X=Y) = \frac{\mu_{12}}{\mu_1+\mu_2+\mu_{12}};
\end{equation}

    \item[(5)] Based on the survival function $\bar{H}$, it is easy to get that the minimum of $X$ and $Y$ is exponential with rate $\mu_1+\mu_2+\mu_{12}$:
\begin{equation}\label{eqn:min}
  \min (X,Y) \sim \mathrm{Exp}(\mu_1+\mu_2+\mu_{12});
\end{equation}

    \item[(6)]  The correlation coefficient $\rho_{X,Y}$ between $X$ and $Y$ is given by
\[
    \rho_{X,Y} = \frac{\mu_{12}}{\mu_1+\mu_2+\mu_{12}}.
\]
It is worthwhile to note that the MO-BVE r.v.s can only be a model for service times with positive dependence;

    \item[(7)] The conditional density $f_{Y|X}(y|x)$ of $Y$ given $X$ is given by
\begin{equation}\label{eqn:density-Y|X}
  f_{Y|X}(y|x) = \left \{ \begin{array}{ll}
    \frac{\mu_1(\mu_2+\mu_{12})}{\mu_1+\mu_{12}} e^{-(\mu_2+\mu_{12})y + \mu_{12}x}, & \mbox{for } y>x, \\
    \mu_2 e^{-\mu_2 x}, & \mbox{for } y<x. \end{array} \right.
\end{equation}
Symmetrically, we have
\begin{equation}\label{eqn:density-X|Y}
  f_{X|Y}(x|y) = \left \{ \begin{array}{ll}
    \frac{\mu_2(\mu_1+\mu_{12})}{\mu_2+\mu_{12}} e^{-(\mu_1+\mu_{12})x + \mu_{12}y}, & \mbox{for } x>y, \\
    \mu_1 e^{-\mu_1 x}, & \mbox{for } x<y. \end{array} \right.
\end{equation}
See equation (47.52) in \cite{Kotz-et at:2000} for more information;

\item[(8)] $P(X <Y)= \mu_1/(\mu_1+\mu_2+\mu_{12})$.
\end{description}
\end{property}

\begin{remark}
It is worthwhile to mention that among all properties of MO-BVED, the property in (5) is a crucial one in the of Theorem~\ref{the:1}.
\end{remark}

\subsection{Construction of MO-BVED via a shock model}

The MO-BVED can be constructed through a shock model, and this construction is important for describing the sample path construction of the $M/M_D/2$ model defined in the next section. 

Consider two components $C_1$ (Type--A) and $C_2$ (Type--B) with three independent
shock clocks started at time $0$: 
$T_1\sim\mathrm{Exp}(\mu_1)$ (kills only $C_1$), 
$T_2\sim\mathrm{Exp}(\mu_2)$ (kills only $C_2$), and
$T_{12}\sim\mathrm{Exp}(\mu_{12})$ (kills both).
Define the lifetimes of $C_1$ and $C_2$ by
\[
X_1=\min\{T_1,T_{12}\},\qquad X_2=\min\{T_2,T_{12}\}.
\]
Then for $x_1,x_2\ge0$,
\[
\overline F(x_1,x_2)
=\Pr(X_1>x_1,X_2>x_2)
=\mathrm{Exp} \bigl(-\mu_1 x_1-\mu_2 x_2-\mu_{12}\max\{x_1,x_2\}\bigr),
\]
so $(X_1,X_2)$ has the MO--BVED
with parameters $(\mu_1,\mu_2,\mu_{12})$.

\begin{remark}[Diagonal singularity]
Common shocks produce a singular mass on the diagonal:
$\Pr(X_1=X_2)>0$, corresponding to simultaneous failures (departures) when the
first shock among $\{T_1,T_2,T_{12}\}$ is the common one.
\end{remark}

\begin{remark}[Queueing interpretation]
If two services start at the same time by the two servers, $X_i$ is the service time
at server $i$. A common shock corresponds to simultaneous service completions.
\end{remark}

\paragraph{Shock \emph{replacement} process:} For the purpose of describing the sample path construction of the $M/M_d/2$ model, we introduce the shock replacement process. 

We adopt a \emph{replacement-as-new policy}: whenever a component fails, that
component is immediately replaced by a new, independent unit; the other, still-working
component continues unchanged. At a common shock (both fail simultaneously), both
components are replaced by new units. By memorylessness of exponentials, after any
failure the remaining times of the still-running clocks are again exponential with
the same rates $(\mu_1,\mu_2,\mu_{12})$.

Starting at time $0$ with both components new, there are four possible orderings among the shocks:

\paragraph{Case (a): $T_{12}<\min\{T_1,T_2\}$.}
A common shock occurs first; both components fail at $T_{12}$. This is a regenerative
epoch: replace both components and restart independent clocks
$T_1',T_2',T_{12}'$ with the same rates.

\paragraph{Case (b): $T_1<T_{12}<T_2$.}
$C_1$ fails at $T_1$; replace it and start a fresh $T_1'\sim\mathrm{Exp}(\mu_1)$.
$C_2$ continues; the common-shock clock keeps running (its residual is exponential and
independent of the past). Continue (possibly multiple new $T_1$ clocks) until the common shock occurs at $T_{12}$; then both
fail and we regenerate as in case (a).

\paragraph{Case (c): $T_2<T_{12}<T_1$.}
Symmetric to (b) with $1\leftrightarrow2$.

\paragraph{Case (d): $\max\{T_1,T_2\}<T_{12}$.}
There are one or more single-component failures before the common shock. After each
single failure, immediately replace that component and continue. When the common shock
arrives (say at time $T_{12}$), both fail and we regenerate.

\begin{remark}[Between regenerations]
Between two consecutive common shocks, there are no simultaneous failures; there may be
multiple alternating single-component failures, each triggering a fresh exponential clock
for the failed component. At each common shock, both are replaced and the process restarts
independently (regeneration).
\end{remark}

\begin{figure}[!htbp]
\centering
\begin{tikzpicture}[x=1.5cm,y=1cm,>=stealth]
\tikzset{
  seg/.style={thick},
  dseg/.style={thick,dashed},
  shock/.style={red,thick},
  tmark/.style={font=\footnotesize},
  caselabel/.style={anchor=west,font=\small,text width=11cm}
}

\begin{scope}[shift={(0,0)}]
  \node[caselabel] at (0,1.2)
    {Case (a): $T_{12}<T_1,T_2$. A common shock at $T_{12}$ kills both components (regeneration).};

  \node at (-0.3,0.35) {$C_1$};
  \node at (-0.3,-0.15) {$C_2$};

  \def\twelveA{4}
  \draw[seg] (0,0.35) -- (\twelveA,0.35);
  \draw[seg] (0,-0.15) -- (\twelveA,-0.15);

  \draw[shock] (\twelveA,0.35) -- (\twelveA,-0.15);
  \node at (\twelveA,0.55) {\(\times\)};
  \node at (\twelveA,-0.35) {\(\times\)};

  \draw[->,seg] (0,-0.75) -- (6,-0.75) node[anchor=west] {time};
  \node[tmark] at (\twelveA,-0.95) {$T_{12}$};
  \node[tmark] at (\twelveA+0.7,-0.95) {$T_1$};
  \node[tmark] at (\twelveA+1.2,-0.95) {$T_2$};
\end{scope}

\begin{scope}[shift={(0,-3.6)}]
  \node[caselabel] at (0,1.4)
    {Case (b): $T_1<T_{12}<T_2$. Type--1 shock at $T_1$ kills $C_1$; $C_1$ is replaced (dashed).
     The common shock at $T_{12}$ kills both (regeneration).};

  \node at (-0.3,0.35) {$C_1$};
  \node at (-0.3,-0.15) {$C_2$};

  \def\toneB{1.2}
  \def\twelveB{4}

  \draw[seg]  (0,0.35) -- (\toneB,0.35);
  \node at (\toneB,0.55) {\(\times\)};
  \draw[dseg] (\toneB,0.35) -- (\twelveB,0.35);

  \draw[seg] (0,-0.15) -- (\twelveB,-0.15);

  \draw[shock] (\twelveB,0.35) -- (\twelveB,-0.15);
  \node at (\twelveB,0.55) {\(\times\)};
  \node at (\twelveB,-0.35) {\(\times\)};
  \node at (\twelveB - 0.5,0.55) {\(\times\)};
  \node at (\twelveB - 1.8,0.55) {\(\times\)};

  \draw[->,seg] (0,-0.75) -- (6,-0.75);
  \node[tmark] at (\toneB,-0.95) {$T_1$};
  \node[tmark] at (\twelveB,-0.95) {$T_{12}$};
  \node[tmark] at (\twelveB+0.9,-0.95) {$T_2$};
\end{scope}

\begin{scope}[shift={(0,-7.2)}]
  \node[caselabel] at (0,1.4)
    {Case (c): $T_2<T_{12}<T_1$. Type--2 shock at $T_2$ kills $C_2$; $C_2$ is replaced (dashed).
     The common shock at $T_{12}$ kills both (regeneration).};

  \node at (-0.3,0.35) {$C_1$};
  \node at (-0.3,-0.15) {$C_2$};

  \def\ttwoC{0.8}
  \def\twelveC{4}

  \draw[seg] (0,0.35) -- (\twelveC,0.35);

  \draw[seg]  (0,-0.15) -- (\ttwoC,-0.15);
  \node at (\ttwoC,-0.35) {\(\times\)};
  \draw[dseg] (\ttwoC,-0.15) -- (\twelveC,-0.15);

  \draw[shock] (\twelveC,0.35) -- (\twelveC,-0.15);
  \node at (\twelveC,0.55) {\(\times\)};
  \node at (\twelveC,-0.35) {\(\times\)};
  \node at (\twelveC - 0.5,-0.35) {\(\times\)};
  \node at (\twelveC - 2.2,-0.35) {\(\times\)};
  \node at (\twelveC - 2.4,-0.35) {\(\times\)};

  \draw[->,seg] (0,-0.75) -- (6,-0.75);
  \node[tmark] at (\ttwoC,-0.95) {$T_2$};
  \node[tmark] at (\twelveC,-0.95) {$T_{12}$};
  \node[tmark] at (\twelveC+0.9,-0.95) {$T_1$};
\end{scope}

\begin{scope}[shift={(0,-11.2)}]
  \node[caselabel] at (0,1.6)
    {Case (d): $T_1<T_2<T_{12}$. Single-component failures occur before the common shock.
     Each failed component is replaced (dashed). At $T_{12}$ both fail (regeneration).};

  \node at (-0.3,0.35) {$C_1$};
  \node at (-0.3,-0.15) {$C_2$};

  \def\toneD{1.2}
  \def\ttwoD{1.7}
  \def\twelveD{4}

  \draw[seg]  (0,0.35) -- (\toneD,0.35);
  \node at (\toneD,0.55) {\(\times\)};
  \node at (\twelveD-0.3,-0.35) {\(\times\)};
  \node at (\twelveD-0.7,-0.35) {\(\times\)};
  \node at (\twelveD-1.2,-0.35) {\(\times\)};
  \draw[dseg] (\toneD,0.35) -- (\twelveD,0.35);

  \draw[seg]  (0,-0.15) -- (\ttwoD,-0.15);
  \node at (\ttwoD,-0.35) {\(\times\)};
  \node at (\twelveD,-0.35) {\(\times\)};
  \draw[dseg] (\ttwoD,-0.15) -- (\twelveD,-0.15);

  \draw[shock] (\twelveD,0.35) -- (\twelveD,-0.15);
  \node at (\twelveD,0.55) {\(\times\)};
  \node at (\twelveD-0.5,0.55) {\(\times\)};

  \draw[->,seg] (0,-0.75) -- (6,-0.75);
  \node[tmark] at (\toneD,-0.95) {$T_1$};
  \node[tmark] at (\ttwoD,-0.95) {$T_2$};
  \node[tmark] at (\twelveD,-0.95) {$T_{12}$};
\end{scope}
\end{tikzpicture}

\caption{Failure patterns ($\times$) under different orderings of $T_1,T_2,T_{12}$
in the Marshall--Olkin shock model. Solid segments show original lifetimes; dashed
segments show lifetimes after replacement. Common-shock times are regeneration epochs.}
\label{fig:MO-cases}
\end{figure}

\section{The $M/M_D/2$ model} \label{sec:3}

In this section, we propose a new queueing model, that can be considered a generalization of the $M/M/c$ queue for $c=2$, denoted by $M/M_D/2$, where $M_D$ is used to emphasize that the service times provided by the two servers can be dependent. Specifically, this is a queueing system with a waiting space of infinite capacity, where the arrivals to this system follow a Poisson process, independent of service times, with rate $\lambda$, and the service times provided by the two servers follow the MO-BVED, given by 
(\ref{eqn:density of MO-BED-b}), where $\lambda_1=\mu_1 + \mu_{12}$ and $\lambda_2=\mu_2 + \mu_{12}$ are the marginal exponential service rates of server 1 and server 2, respectively, and $\mu_{12}$ is the dependence parameter. 

Define the state to be the number of customers in the system. The state space 
 $S$ is defined as
\begin{equation}\label{eqn:state-space-MMD2}
    S=\{0, (1,0), (0,1), 2, 3, \ldots \},
\end{equation}
where state $n$, for $n =0$ and $n \geq 2$, represents the number of customers in the system, and states $(1,0)$ and $(0,1)$ represent only one customer in the system, being served by server 1 and 2, respectively.

\begin{remark}
Since the marginal service rates are different in general, we cannot use 1 to represent the state of one customer in the system. Instead, in this case, we need to specify which server is busy using $(1,0)$ or $(0,1)$. This definition of the state space was also used in modelling the $M/M/2$ queueing system with two heterogeneous servers (for example, see \cite{Gumbel:1960} and \cite{Singh:1970}). 
\end{remark}

\begin{remark}
Intuitively, the $M/M_D/2$ model is stable if and only if $\lambda < \mu_1 + \mu_2 + \mu_{12}$, which can be proved through identifying a Lyapunov function using the standard drift condition for Markov chains. 
This condition is intuitively clear since stability requires that the arrival rate ($\lambda$) be smaller than the maximum average departure rate when both servers are busy. The departure rate consists of two components: (1) the two servers working individually, removing customers at total rate $\mu_1+\mu_2$, and (2) a bonus mechanism that clears both customers simultaneously at rate $\mu_{12}$, contributing $2\mu_{12}$ to the departure rate.
\end{remark}

An arrival to the empty system would change the state to $(1,0)$ with probability $0 \leq p \leq 1$, or to $(0,1)$ with probability $q=1-p$; an arrival to state $(1,0)$ or $(0,1)$ changes the state to $2$; and an arrival to state $n$ for $n \geq 2$ changes the state to $n+1$. When the system is not empty, we have the following two cases:

\paragraph{Case~I: One busy server.}
Suppose server $i\in\{1,2\}$ is busy and the other is idle. Generate
$T_i\sim\mathrm{Exp}(\mu_i)$ and $T_{12}\sim\mathrm{Exp}(\mu_{12})$, independent; the
service time is $X_i=\min\{T_i,T_{12}\}$. If the next arrival occurs before the service completion of server $i$,
it begins service immediately at server $\bar i$ (the server other than $i$) with its own clock
$T_{\bar i}\sim\mathrm{Exp}(\mu_{\bar i})$ and its service time is $X_{\bar{i}} = \min\{T_{\bar{i}}, \text{ remaining time of }T_{12}\}$; from that instant on, the remaining service
times (note that a new service time can be treated as a remaining service time) of the two busy servers form an MO--BVED with parameters
$(\mu_1,\mu_2,\mu_{12})$.

\paragraph{Case~II: Two busy servers (with or without a queue).}
Whenever both servers are busy, the vector of remaining service times has MO--BVE law
with parameters $(\mu_1,\mu_2,\mu_{12})$. A type-$i$ (single) shock completes only
server $i$; a common shock completes both simultaneously. Arrivals, if any, join the
queue. After each completion, memorylessness ensures that if both servers are still
busy, their new remaining times are again MO--BVE.

\begin{remark}
When a customer arrives to the empty system, its service time is $X_i=\min\{T_i,T_{12}\}$. Although the distribution of $X_i$ is $\mathrm{Exp}(\mu_i + \mu_{12})$ (the marginal distribution of server $i$), the marginal distribution does not provide detailed information about the clock time $T_{12}$, which is needed for determining the service time at server $\bar{i}$ if the next arrival is before the departure from server $i$. 
\end{remark}

\begin{remark}[Coupling to the shock replacement process]
For $n\ge4$ customers in system (both busy plus a queue), the busy-period segment between
consecutive common shocks follows exactly the shock replacement dynamics described above:
single-completion events alternate with occasional simultaneous completions at common
shocks (regenerations).
\end{remark}

The above description provides details of the connection of the model to MO-BVED, and therefore a sample path construction of the model. We will show in the next subsection that  $\{ X(t): t \geq 1\}$ is a Markov chain.

It is our expectation that this new queueing system could be used as a mathematical model for applications, in which service times are positively correlated. We hope that this basic system can play a role in modelling queueing systems with dependent servers, similar to the role played by the $M/M/c$ queueing for systems with independent servers.

\subsection{Markov chain}

In this section, we prove that the process $X(t)$ defined in the previous section with state space $S$ specified in (\ref{eqn:state-space-MMD2}) is a Markov chain. It is worth noting that modelling a queueing system with dependent service times as a Markov chain is a new approach proposed in this paper. With this Markov chain, we hope many classical methods applied to queueing systems with independent service times can be extended to queueing systems with dependent service times.

\begin{theorem} \label{the:1}
For the $M/M_D/2$ queueing system with arrival rate $\lambda$ and two heterogeneous dependent servers characterized by the MO-BVED with parameters $\mu_1$, $\mu_2$ and $\mu_{12}$ given in (\ref{eqn:density of MO-BED-b}), the process   $\{X(t): t \geq 0\}$ defined in the previous section, representing the number of customers in the system at time $t$, is a continuous time Markov chain with transition rate matrix $Q = (q_{i_1,i_2})$ given by:
\begin{equation}\label{eqn:Q-MMD2}
    Q = \bordermatrix[{[]}]{ & 0 & (1,0) & (0,1) &2 &3 & 4 & 5 & 6 & \cdots & &\cr
    0 & -\lambda & p\lambda & q\lambda & \cr
    \hspace*{-3mm}(1,0) & \mu_1+\mu_{12} & -\theta_{(1,0)} & 0 & \lambda &\cr
    \hspace*{-3mm}(0,1) & \mu_2+\mu_{12} & 0 & -\theta_{(0,1)} & \lambda &\cr
    2 & \mu_{12} & \mu_2 & \mu_1 & - \theta & \lambda & \cr
    3 & & p \mu_{12} & q \mu_{12} & \mu_1+\mu_2 & -\theta & \lambda & \cr
    4 & & & & \mu_{12} & \mu_1+\mu_2 & -\theta & \lambda & \cr
    5 & & & & & \mu_{12} & \mu_1+\mu_2 & -\theta & \lambda & \cr
    \vdots & & &&&&\ddots & \ddots & \ddots & \ddots },
\end{equation}
where all empty entries are 0, $\theta_{(1,0)}=\lambda+\mu_1+\mu_{12}$, $\theta_{(0,1)}=\lambda+\mu_2+\mu_{12}$,  $\theta=\lambda+\mu_1+\mu_2+\mu_{12}$, and $0 < p <1$ with $q=1-p$. 

\end{theorem}

\proof We divide the proof into two steps: in the first step, we show that the sojourn time in each state is exponential with the parameter dependent on the current state (independent of the future state); and in the second step, we explicitly derive the transition probabilities for the imbedded Markov chain.

\textbf{Step~1.} 

\textbf{State 0:} For $i=0$,
\[
    \mbox{sojourn time in 0} = \mbox{remaining interarrival time} = A \sim \mathrm{Exp}(\theta_0),
\]
where $A$ is the interarrival time and $\theta_0=\lambda$.

\textbf{State (1,0):} For state $(1,0)$,
\begin{align*}
    \mbox{sojourn time in $(1,0)$} & = \min (\mbox{remaining interarrival time, remaining service time at server 1)} \\
     &= \min (A, S_1) \sim \mathrm{Exp}(\theta_{(1,0)}),
\end{align*}
where $\theta_{(1,0)} = \lambda + \mu_1+\mu_{12}$ since only server~1 is busy and the marginal service time $S_1 \sim \mathrm{Exp}(\mu_1+\mu_{12})$.

\textbf{State (0,1):} For state $(0,1)$, it can be similarly shown that
\begin{align*}
    \mbox{sojourn time in $(0,1)$} & = \min (\mbox{remaining interarrival time, remaining service time at server 2)} \\
     &= \min (A, S_2) \sim \mathrm{Exp}(\theta_{(0,1)}),
\end{align*}
where $\theta_{(0,1)} = \lambda + \mu_2+\mu_{12}$ since only server~2 is busy and the marginal service time $S_2 \sim \mathrm{Exp}(\mu_2+\mu_{12})$.

\textbf{State $i$:} For state $i \geq 2$,
\begin{align*}
    \mbox{sojourn time in $i$} &= \min (\mbox{remaining interarrival time, remaining service times at servers 1 and 2)} \\
    &= \min (A, (S_1, S_2)) = \min(A, \min(S_1,S_2)) \sim \mathrm{Exp}(\theta_i),
\end{align*}
where $\theta_i = \theta = \lambda + \mu_1+\mu_2 + \mu_{12}$ since $\min(S_1,S_2) \sim \mathrm{Exp}(\mu_1+\mu_2+\mu_{12})$.
\medskip

\textbf{Step~2.}
Let $X_n$  be the imbedded process of $X(t)$, or $X_n$ is the state to which $X(t)$ moves when it makes its $n$th transition. For the model with heterogeneous servers, when a customer is arriving to an empty system, we need to specify which server the arriving customer would be sent to. We assume that with probability $p$, the process moves from 0 to $(1,0)$, and with $q=1-p$ from 0 to $(0,1)$.

\textbf{State 0:} When $X_n=0$, it follows from the above assumption that
\[
    p_{0,(1,0)} = p \quad \mbox{and} \quad p_{0,(0,1)} = q =1-p,
\]
which are independent of past transitions.

\textbf{State (1,0):} When $X_n=(1,0)$,
\begin{align*}
  P(X_{n+1} = 0 \vert X_n =(1,0)) &= P(\mbox{remaining service time at server 1 $<$ remaining interarrival time}) \\
 &= P(S_1<A) = \frac{\mu_1+\mu_{12}}{\lambda+\mu_1+\mu_{12}} = p_{(1,0),0},
\end{align*}
which is independent of past transitions; and similarly
\begin{align*}
 P(X_{n+1} = 2 \vert X_n =(1,0)) &= P(\mbox{remaining service time at server 1 $>$ remaining interarrival time}) \\
 &=P(A<S_1) = \frac{\lambda}{\lambda+\mu_1+\mu_{12}} = p_{(1,0),2}.
\end{align*}

\textbf{State (0,1):} When $X_n = (0,1)$, we can similarly show that
\[
    p_{(0,1),0} = \frac{\mu_2+\mu_{12}}{\lambda+\mu_2+\mu_{12}} \quad \mbox{and} \quad p_{(0,1),2} = \frac{\lambda}{\lambda+\mu_2+\mu_{12}}.
\]

\textbf{State 2:} When $X_n = 2$,
\begin{align*}
 P(X_{n+1} = 3 \vert X_n =2) &= P(A < \min(S_1, S_2)) \\
 &=  \frac{\lambda}{\lambda+\mu_1+\mu_2+\mu_{12}} = p_{2,3},
\end{align*}
which is independent of past transitions.

Note that the MO-BVED is singular on $S_1=S_2$, which implies that with a positive probability we can have simultaneous service completions from both servers.
Noticing that
\[
    P(S_1=S_2 \leq t) = \frac{\mu_{12}}{\mu_1+\mu_2+\mu_{12}} \big [1 - e^{-(\mu_1+\mu_2+\mu_{12})t} \big ], \quad t \geq 0,
\]
according to (\ref{eqn:X=Y>t}) and (\ref{eqn:X=Y}), which leads to the density $f_{S_1=S_2}(t)$ of $S_1=S_2$, given by
\begin{equation}\label{eqn:density-X=Y}
    f_{S_1=S_2}(t) = \mu_{12}\, e^{-(\mu_1+\mu_2+\mu_{12})t}, \quad t \geq 0.
\end{equation}
Then, by conditioning on the value of $S_1=S_2$, we have
\begin{align*}
    p_{2,0} &= P(S_1=S_2<A) = \int_{0}^{\infty} \mu_{12}\, e^{-(\mu_1+\mu_2+\mu_{12})t} P(A > t \vert S_1=S_2=t) dt \\
    &= \int_{0}^{\infty} \mu_{12}\, e^{-(\mu_1+\mu_2+\mu_{12})t} P(A > t) dt \\
    &= \mu_{12} \int_{0}^{\infty} e^{-(\mu_1+\mu_2+\mu_{12})t} e^{-\lambda t} dt \\
    &= \frac{\mu_{12}}{\lambda+\mu_1+\mu_2+\mu_{12}} \int_{0}^{\infty} (\lambda+\mu_1+\mu_2+\mu_{12}) e^{-(\lambda+\mu_1+\mu_2+\mu_{12})t} dt \\
    & = \frac{\mu_{12}}{\lambda+\mu_1+\mu_2+\mu_{12}}.
\end{align*}

By conditioning on the value of $S_2$, we now calculte
\begin{align*}
 p_{2,(1,0)}&= P(\min(S_1,S_2) <A, S_2 < S_1) = P(S_2 <A, S_2 < S_1)\\
 &= \int_{0}^{\infty} f_{S_2}(t) P(S_2 <A, S_2 < S_1 \vert S_2=t) dt \\
  &= \int_{0}^{\infty} f_{S_2}(t) P(A >t, S_1 >t \vert S_2=t) dt \\
   &= \int_{0}^{\infty} f_{S_2}(t) P(A >t\vert S_2=t) P(S_1 >t \vert A>t, S_2=t) dt \\
    &= \int_{0}^{\infty} f_{S_2}(t) P(A >t) P(S_1 >t \vert S_2=t) dt.
\end{align*}
Notice that the density function $f_{S_1|S_2}(x|y)$ is given by (\ref{eqn:density-X|Y}), or
\[
  f_{X|Y}(x|y) = \frac{\mu_2(\mu_1+\mu_{12})}{\mu_2+\mu_{12}} e^{-(\mu_1+\mu_{12})x + \mu_{12}y}, \quad \mbox{for } x>y,
\]
based on which we calculate
\begin{align*}
    P(S_1 >t \vert S_2=t) &= \int_{t}^{\infty} \frac{\mu_2(\mu_1+\mu_{12})}{\mu_2+\mu_{12}} e^{-(\mu_1+\mu_{12})x + \mu_{12}t} dx \\
     & =  \frac{\mu_2(\mu_1+\mu_{12})}{\mu_2+\mu_{12}} e^{\mu_{12} t} \int_{t}^{\infty} e^{-(\mu_1+\mu_{12})x} dx \\
     &= \frac{\mu_2}{\mu_2+\mu_{12}} e^{\mu_{12} t} e^{-(\mu_1+\mu_{12})t} = \frac{\mu_2}{\mu_2+\mu_{12}} e^{-\mu_{1}t}.
\end{align*}

Finally, by noticing that $f_{S_2}(t)$ is exponential with parameter $\mu_2+\mu_{12}$, we have
\begin{align*}
     p_{2,(1,0)}&= \int_{0}^{\infty} (\mu_2+\mu_{12}) e^{-(\mu_2+\mu_{12})t} \, \cdot \, e^{-\lambda t} \, \cdot \, \frac{\mu_2}{\mu_2+\mu_{12}} e^{-\mu_{1}t} dt \\
     &= \frac{\mu_{2}}{\lambda+\mu_1+\mu_2+\mu_{12}}.
\end{align*}

Similarly,
\begin{align*}
 P(X_{n+1} = (0,1) \vert X_n =2) & = P(\min(S_1,S_2) <A, S_1 < S_2) = P(S_1 <A, S_1 < S_2)\\
 &= \frac{\mu_{1}}{\lambda+\mu_1+\mu_2+\mu_{12}}.
\end{align*}

\textbf{State $i$:} When $X_n =i \geq 3$,
\begin{align*}
 P(X_{n+1} = i+1 \vert X_n =i) &= P(A < \min(S_1, S_2)) \\
 &=  \frac{\lambda}{\lambda+\mu_1+\mu_2+\mu_{12}} = p_{i,i+1},
\end{align*}
which is independent of past transitions.

Similarly,
\begin{align*}
 P(X_{n+1} = i-1 \vert X_n =i) &=  P(A > \min(S_1, S_2), S_1 \neq S_2) \\
 &=  \frac{\mu_1+\mu_2}{\lambda+\mu_1+\mu_2+\mu_{12}} = p_{i,i-1},
\end{align*}
which is independent of past transitions.

We now consider simultaneous transitions. In this case, we need to distinguish $i=3$ from $i>3$. For $i>3$, by conditioning on $S_1=S_2=t$ and using the density $f_{S_1=S_2}(t)$ in (\ref{eqn:density-X=Y}), we have
\begin{align*}
 P(X_{n+1} = i-2 \vert X_n =i) &=  P(S_1=S_2 < A) \\
 &   =\int_{0}^{\infty} \mu_{12} e^{-(\mu_1+\mu_2+\mu_{12})t} e^{-\lambda t} dt,\\
 &= \frac{\mu_{12}}{\lambda+\mu_1+\mu_2+\mu_{12}} = p_{i,i-2}.
\end{align*}

For $i=3$, if there are two simultaneous service completions, then the number of customers in the system is reduced to 1. In this case, with probability $p$ and $q=1-p$, the Markov chain enters state $(1,0)$ and $(0,1)$, respectively, or,
\[
    P(X_{n+1} = (1,0) \vert X_n =3) = \frac{p \mu_{12}}{\lambda+\mu_1+\mu_2+\mu_{12}} = p_{3,(1,0)},
\]
and
\[
    P(X_{n+1} = (0,1) \vert X_n =3) = \frac{q \mu_{12}}{\lambda+\mu_1+\mu_2+\mu_{12}} = p_{3,(0,1)}.
\]

The transition rate matrix follows now immediately from 
$q_{i_1,i_1} = - \theta_{i_1}$ and $q_{i_1,i_2} = \theta_{i_1} p_{i_1,i_2}$ for $i_1 \neq i_2$.
\pend


\begin{remark} 
\textbf{(1)} For the $M/M_D/2$ model with two heterogeneous servers, we cannot combine states $(1,0)$ and $(0,1)$ into a single state 1. Otherwise, the sojourn time in state 1 is not exponential. For the case of homogeneous servers, we can introduce state 1 (see next section). 

\textbf{(2)} The continuous-time Markov chain in~(\ref{eqn:Q-MMD2}) is not a birth-and-death process. This is because simultaneous departures from both servers are possible, a consequence of singularities in the MO-BVE service time density along the diagonal. This is a novel phenomenon that does not arise in the Markov chain for the $M/M/c$ model.

\end{remark}

\subsection{$M/M_D/2$ model with homogeneous servers}

The case of $\mu_1=\mu_2=\mu$ is special. In this case, we can combine states $(1,0)$ and $(0,1)$ into a single state 1. We now modify the state space $S$ to define 
\begin{equation}\label{eqn:S0}
   S_0 =\{0,1,2, \ldots\}.
\end{equation}
Then, the process $X(t)$, defined on the space $S_0$, representing the number of customers in the system at time $t$, is a Markov chain, as stated in the following corollary.

\begin{corollary} \label{cor:1}
The process $X(t)$ defined on $S_0$ for the $M/M_D/2$ queueing system with dependent service times characterized by the MO-BVED with $\mu_1=\mu_2$, denoted by $\mu$, is a continuous time Markov chain, and the transition rate matrix $Q_0 = (q_{ij})$ with $q_{ii} = - \theta_i$ and $q_{ij} = \theta_i p_{ij}$ for $i \neq j$ is given by
\begin{equation}\label{rate-matrix-Q0}
        Q_0 = \bordermatrix[{[]}]{ & 0 & 1 &2 &3 & 4 & 5 & \cdots & \cr
    0 & -\lambda & \lambda & 0 \cr
    1 & \mu+\mu_{12} & -(\lambda+\mu+\mu_{12}) & \lambda \cr
    2 & \mu_{12} & 2\mu & -\theta_2 & \lambda \cr
    3 & & \mu_{12} & 2\mu & -\theta_3 & \lambda \cr
    4 & & & \mu_{12} & 2\mu & -\theta_3 & \lambda \cr
     \vdots & & && \ddots & \ddots & \ddots & \ddots },
\end{equation}
where, for $i \geq 2$, $\theta_i = \lambda+2\mu+\mu_{12}$.
\end{corollary}

\proof The proof for case of heterogeneous servers is still valid with the following modifications:

When $i=1$,
\[
    \mbox{the sojourn time in 1 $\sim \mathrm{Exp}(\theta_1)$ with $\theta_1=\lambda + (\mu+\mu_{12})$}
\]
by noticing that the marginal service time for each server is $\mathrm{Exp}(\mu+\mu_{12})$. It is clear that
$p_{01} = 1$.

For $X_n=1$, since the marginal distributions are the same, which is $\mathrm{Exp}(\mu+\mu_{12})$, we have
\begin{align*}
 p_{10}= P(X_{n+1} = 0 \vert X_n =1) &= P(\mbox{remaining service time $<$ remaining interarrival time}) \\
 &=P(S_1 \mbox{ (or $S_2$)} <A) = \frac{\mu+\mu_{12}}{\lambda+\mu+\mu_{12}},
\end{align*}
which is independent of past transitions;

Similarly,
\begin{align*}
 p_{12}= P(X_{n+1} = 2 \vert X_n =1) &= P(\mbox{remaining service time $>$ remaining interarrival time}) \\
 &=P(A<S_1 \mbox{ (or $S_2$)}) = \frac{\lambda}{\lambda+\mu+\mu_{12}},
\end{align*}
\[
    p_{21} = p_{32} = 2\mu/(\lambda+2\mu+\mu_{12})
\]
and
\[
    p_{31}= p_{20} = \frac{\mu_{12}}{\lambda+2\mu+\mu_{12}}.
\]

The expression of the transition rate matrix $Q_0$ follows immediately from the above results. \pend

\section{Stationary distribution of the $M/M_D/2$ queue} \label{sec:4}


Assume that the Markov chain has a unique stationary probability vector $\pi=(\pi_0, \pi_{(1,0)}, \pi_{(0,1)}, \pi_2, \pi_3, \ldots, )$ (a necessary and sufficient condition for the system to be stable will be provided later in this section). The purpose is to find the solution of $\pi$ of the stationary equations, given by
\begin{equation}\label{piQ=0}
    \pi Q=0,
\end{equation}
subject to the normalization condition 
\begin{equation}\label{normalization-Q}
    \pi_0 + \pi_{(1,0)} + \pi_{(0,1)} + \sum_{i=2}^\infty \pi_i =1.
\end{equation}
For convenience, we write out in detail all stationary equations:
\begin{align}
  0  = & -\lambda \pi_0 + (\mu_1+\mu_{12}) \pi_{(1,0)} + (\mu_2+\mu_{12}) \pi_{(0,1)} + \mu_{12} \pi_2, \label{eqn-0-for-Q} \\
  0  = & p\lambda \pi_0 - (\lambda + \mu_1 + \mu_{12}) \pi_{(1,0)} +  \mu_2 \pi_2 + p\mu_{12} \pi_3, \label{eqn-(1,0)-for-Q} \\
  0  = & q\lambda \pi_0 - (\lambda + \mu_2 + \mu_{12}) \pi_{(0,1)} +  \mu_1 \pi_2 + q\mu_{12} \pi_3, \label{eqn-(0,1)-for-Q} \\
  0  = & \lambda (\pi_{(1,0)} + \pi_{(0,1)}) - (\lambda + \mu_1 + \mu_2 + \mu_{12}) \pi_2 + (\mu_1+\mu_2) \pi_{3} + \mu_{12} \pi_{4}, \label{eqn-2-for-Q} \\
  0  = & \lambda \pi_{i-1} - (\lambda + \mu_1 + \mu_2 + \mu_{12}) \pi_i +  (\mu_1 + \mu_2) \pi_{i+1} + \mu_{12} \pi_{i+2}, \quad i \geq 3. \label{eqn->2-for-Q}
\end{align}

We use the following standard steps to show that the stationary distribution $\pi$ has a geometric tail:

\paragraph{Step 1: Geometric tail and characteristic equation:}
For geometric tail, assume $\pi_i = C r^i$ for $i \geq 2$, with which the non-boundary equation \eqref{eqn->2-for-Q} leads to the following characteristic equation:
\begin{equation}\label{ch-eqn-Q}
    \mu_{12} r^3 + (\mu_1+\mu_2)  r^2 - (\lambda + \mu_1 + \mu_2 + \mu_{12}) r + \lambda = 0.
\end{equation}
Notice that $r=1$ is a solution. Hence, we can write the characteristic equation as
\begin{equation}
    (r-1)[\mu_{12} r^2 + (\mu_1+\mu_2 + \mu_{12}) r - \lambda] = 0.
\end{equation}
Solving the quadratic equation 
\[
    \mu_{12} r^2 + (\mu_1+\mu_2 + \mu_{12}) r - \lambda =0
\]
gives the two roots:
\begin{align}
    r_- &= \frac{-(\mu_1 + \mu_2 + \mu_{12}) - \sqrt{\Delta}}{2 \mu_{12}}, \\
    r_+ &= \frac{-(\mu_1 + \mu_2 + \mu_{12}) + \sqrt{\Delta}}{2 \mu_{12}}, \label{eqn:r+}
\end{align}
where $\Delta = (\mu_1 + \mu_2 + \mu_{12})^2 + 4 \lambda \mu_{12}$.
It is clear that $r_-<0$ and $|r_-| > 1$. Therefore, for a convergent solution, only $r_+$ (which satisfies $r_+ >0$) can be used. 

\begin{remark} Note that $r_+ < 1$ if and only if  $\lambda < \mu_1 + \mu_2 + 2\mu_{12}$ through simple algebra.
\end{remark}

\paragraph{Step 2: Consistency of geometric tail with non-boundary equations:}
In this step, let $r=r_+$. We show that replacing the geometric tail $\pi_i = Cr^i$ ($i \geq 2$) in the boundary equations \eqref{eqn-0-for-Q} -- \eqref{eqn->2-for-Q} gives unique non-negative determination of the boundary probabilities $\pi_0$, $\pi_{(1,0)}$, $\pi_{(0,1)}$, all of which are functions of $C$. Specifically,

From the equation for state 2, or \eqref{eqn-2-for-Q}, we also have
\[
    \pi_{(1,0)} + \pi_{(0,1)} = C r,
\]  
which implies that
\[
    \pi_{(1,0)} = C r - \pi_{(0,1)}.
\]  
Substitute the above solutions into equation \eqref{eqn-0-for-Q} and equation \eqref{eqn-(0,1)-for-Q} to get:
\[
    \pi_{(0,1)} = C r \frac{A}{B},
\]
where
\begin{align}
    A & = \lambda - p \mu_1 + \mu_1 - p \mu_{12} (r^2 + r + 1) + \mu_{12} - \mu_2 r, \label{eqn:A} \\
    B & = \lambda - p \mu_1 + \mu_1 + \mu_{12} + p \mu_2, \label{eqn:B}
\end{align} 
and
\[
    \pi_0 = \frac{C r}{\lambda B} P(r),
\]
where
\begin{align} \nonumber
    P(r) = & \lambda \mu_{12} r + \lambda \mu_{12} + \lambda \mu_2 + \mu_1 \mu_{12} p r^2 + \mu_1 \mu_{12} r \\ &
     + \mu_1 \mu_{12} + \mu_1 \mu_2 r + \mu_1 \mu_2 + \mu_{12}^2 r + \mu_{12}^2 - \mu_{12} \mu_2 p r^2 + \mu_{12} \mu_2 - \mu_2^2 r. \label{eqn:P(r)}
\end{align}
We then also have
\[
    \pi_{(1,0)} = C r - \pi_{(0,1)} = C r \frac{B-A}{B}.
\]

\begin{remark}
With the above $\pi_0$, $\pi_{(0,1)}$, $\pi_{(1,0)}$ and $\pi_i = C r^i$ for $i \geq 2$, we can check this set of solutions also satisfies \eqref{eqn-(1,0)-for-Q} as expected. In fact, 
The four terms in equation \eqref{eqn-(1,0)-for-Q} become:
\begin{align*}
\text{Term 1:} & \quad p \lambda \pi_0 = \frac{C r \cdot p}{B} \cdot P(r); \\
\text{Term 2:} & \quad -(\lambda + \mu_1 + \mu_{12}) \cdot \pi_{(1,0)} = -(\lambda + \mu_1 + \mu_{12}) C r \cdot \frac{B - A}{B}; \\
\text{Term 3:} & \quad \mu_2 \pi_2 = \mu_2 C r^2; \\
\text{Term 4:} & \quad p \mu_{12} \pi_3 = p \mu_{12} C r^3.
\end{align*}
Combining the terms and simplifying, we find:
\[
    p \lambda \pi_0 - (\lambda + \mu_1 + \mu_{12}) \pi_{(1,0)} + \mu_2 \pi_2 + p \mu_{12} \pi_3 = 0,
\]
which is equation \eqref{eqn-(1,0)-for-Q}.
\end{remark}

After the consistency check, we have the following condition for the stability of the system.
\begin{corollary}
The $M/M_D/2$ system has a unique stationary distribution, if and only if $\lambda < \mu_1 + \mu_2  + 2\mu_{12}$, or $\rho < 1$ with $\rho = \lambda / (\mu_1 + \mu_2  + 2\mu_{12})$.
\end{corollary}

\paragraph{Step 3: Probability solution of $\pi$:}
First, we show that $\pi_0$, $\pi_{(0,1)}$ and $\pi_{(1,0)}$ are all positive. Note that
\[
B=\lambda+(1-p)\mu_1+\mu_{12}+p\mu_2>0.
\]
Moreover,
\[
B-A=p\mu_2+\mu_2 r+p\mu_{12}(r^2+r+1)>0.
\]
Hence
\[
\pi_{(1,0)}=Cr\frac{B-A}{B} > 0.
\]

Next, using the root relation
\[
\lambda=\mu_{12}r^2+(\mu_1+\mu_2+\mu_{12})r,
\]
we obtain
\[
A
=\mu_1 r+(1-p)\bigl(\mu_{12}r^2+\mu_{12}r+\mu_1+\mu_{12}\bigr)>0.
\]
Therefore
\[
\pi_{(0,1)}=Cr\frac{A}{B} > 0.
\]

Finally, from the balance equation at state \(0\),
\[
\lambda\pi_0=(\mu_1+\mu_{12})\pi_{(1,0)}+(\mu_2+\mu_{12})\pi_{(0,1)}+\mu_{12}\pi_2,
\]
and hence
\[
\pi_0
=
\frac{Cr}{\lambda B}
\Bigl[(\mu_1+\mu_{12})(B-A)+(\mu_2+\mu_{12})A+\mu_{12}rB\Bigr] > 0,
\]
when $C >0$.

To get the probability solution, we use the normalization condition \eqref{normalization-Q}, which is equivalent to 
\[
   \frac{C r}{\lambda B} P(r) + C r + \frac{C r^2}{1 - r} = 1,
\]
to determine the value of $C$. We summarize the above results into the following theorem.
\begin{theorem} \label{the:2}
For the $M/M_D/2$ queueing system satisfying the stability condition $\lambda < \mu_1+\mu_2 + 2\mu_{12}$, the stationary probability vector is given by
\begin{align}
    \pi_0 & = \frac{C r}{\lambda B} P(r), \\
    \pi_{(0,1)} & = C r \frac{A}{B}, \\
    \pi_{(1,0)} & =  C r \frac{B-A}{B}, \\
    \pi_k & = C r^k, \quad k \geq 2,
\end{align}
where 
\begin{align}
    A = & \lambda - p \mu_1 + \mu_1 - p \mu_{12} (r^2 + r + 1) + \mu_{12} - \mu_2 r, \\
    B = & \lambda - p \mu_1 + \mu_1 + \mu_{12} + p \mu_2, \\
        P(r)  = & \lambda \mu_{12} r + \lambda \mu_{12} + \lambda \mu_2 + \mu_1 \mu_{12} p r^2 + \mu_1 \mu_{12} r \\ &
     + \mu_1 \mu_{12} + \mu_1 \mu_2 r + \mu_1 \mu_2 + \mu_{12}^2 r + \mu_{12}^2 - \mu_{12} \mu_2 p r^2 + \mu_{12} \mu_2 - \mu_2^2 r, \\
     C = &  \left [ \frac{r}{\lambda B} P(r) +  \frac{r}{1 - r} \right ]^{-1}, \label{eqn:C}
\end{align} 
and 
\begin{equation}\label{r+-for-Q}
    r = \frac{-(\mu_1 + \mu_2 + \mu_{12}) + \sqrt{\Delta}}{2 \mu_{12}},
\end{equation}
where 
\begin{equation}\label{eqn:Delta}
\Delta = (\mu_1 + \mu_2 + \mu_{12})^2 + 4 \lambda \mu_{12}.
\end{equation}
\end{theorem}

\begin{corollary}
The expected number of customers in the system is given by
\begin{equation}\label{ELD-1}
  E[L_D] = C \cdot \frac{r}{(1-r)^2},
\end{equation}
and the expected number of customers waiting in the queue is given by
\begin{equation}\label{EQD-1}
  E[Q_D] = C \cdot \frac{r^3}{(1-r)^2}.
\end{equation}
\end{corollary}

\subsection{$M/M_D/2$ with homogeneous servers}

When $\mu_1 = \mu_2$, states $(1,0)$ and $(0,1)$ collapse to a single state 1. In this case, the transition rate matrix is given by \eqref{rate-matrix-Q0}. The stationary equations are given by
\begin{align}
  0  = & -\lambda \pi_0 + (\mu+\mu_{12}) \pi_1 + \mu_{12} \pi_2, \label{eqn-equal-rate-1} \\
  0  = & \lambda \pi_0 - (\lambda + \mu + \mu_{12}) \pi_1 + 2 \mu \pi_2 + \mu_{12} \pi_3, \label{eqn-equal-rate-2} \\
  0  = & \lambda \pi_{i-1} - (\lambda + 2\mu + \mu_{12}) \pi_i + 2 \mu \pi_{i+1} + \mu_{12} \pi_{i+2}, \quad i \geq 2. \label{eqn-equal-rate-3}
\end{align}
In this case, the solution of $\pi Q_0=0$ subject to the normalization $\sum_{i=0}^{\infty} \pi_i=1$ is given by the following corollary.

\begin{corollary} \label{cor:2}
For the $M/M_D/2$ queueing system with $\mu_1=\mu_2=\mu$ satisfying the stability condition $\lambda < 2(\mu + \mu_{12})$, the stationary probability vector is given by
\begin{align}
    \pi_0 & = C_h  \frac{(\mu + \mu_{12}) r_h + \mu_{12} r_h^2}{\lambda}, \label{eqn:pi0-mu1=mu2} \\
    \pi_i & = C_h r_h^i, \quad i \geq 1,
\end{align}
where 
\begin{equation}\label{C-h-mu1=mu2}
      C_h = \left [\frac{(\mu + \mu_{12}) r_h + \mu_{12} r_h^2}{\lambda} + \frac{r_h}{1 - r_h}  \right ]^{-1}
\end{equation}
and 
\begin{equation}\label{r+-mu1=mu2}
     r_h = \frac{-2\mu - \mu_{12} + \sqrt{4\lambda \mu_{12} + 4\mu^2 + 4\mu\mu_{12} + \mu_{12}^2}}{2\mu_{12}}.
\end{equation}
\end{corollary}

\begin{remark}
It follows from simple calculations that when $\mu_1=\mu_2$, then $B$ given in \eqref{eqn:B} is simplified to
\[
    B = \lambda + \mu + \mu_{12},
\]
and $P(r)$ given in \eqref{eqn:P(r)} is simplified to
\[
    P(r_h) = (\lambda + \mu + \mu_{12}) [\mu_{12}(1+r_h) + \mu].
\]
Therefore, 
\[
    \frac{P(r_h)}{B} = \mu_{12}(1+r_h) + \mu,
\]
from which it is easy to see $C = C_h$.
\end{remark}

\begin{corollary}
The expected number of customers in the system is given by
\begin{equation}\label{ELD}
  E[L_{D,h}] = C_h \cdot \frac{r_h}{(1-r_h)^2},
\end{equation}
and the expected number of customers waiting in the queue is given by
\begin{equation}\label{EQD}
  E[Q_{D,h}] = C_h \cdot \frac{r_h^3}{(1-r_h)^2}.
\end{equation}

\end{corollary}

\proof Simple calculations lead to a proof. \pend

\section{Impact of service-time dependence on performance} \label{sec:5}

It is of interest to see how dependence would impact on the system performance, which is the main focus of this section. 

\subsection{Solution when $\mu_{12} \to 0$}

Before providing numerical results, we point out that intuitively, when $\mu_{12} \to 0$, the solution of the $M/M_D/2$ model is expected to approach to the solution for the $M/M/2$ model with heterogeneous servers. For the case of $\mu_1=\mu_2$, our solution is expected to approach to the solution for the standard $M/M/2$ system, or with homogeneous servers. We now make a confirmation of this intuition. 

\subsubsection{When $\mu_1=\mu_2 =\mu$}

In this case, there are two homogeneous servers in the system, each with service rate $\mu$. If service times by the two servers are independent, the system becomes the standard $M/M/2$ queue, for which the solution is well-known.

\begin{corollary} \label{cor:to-MM2-mu1=mu2}
For the $M/M_D/2$ model with homogeneous servers (or $\mu_1=\mu_2=\mu$), as $\mu_{12} \to 0$, we have the following:
\begin{align}\label{eqn:mu{12}-to-0-when-mu1=mu2}
  r_{0,h} = & \lim_{\mu_{12} \to 0} r_h = \frac{\lambda}{2\mu}, \\
  C_{0,h} = & \lim_{\mu_{12} \to 0} C_h = 2 \frac{1-\rho_{0,h}}{1+\rho_{0,h}} = 2 \pi_0, \\
\end{align}
where $\rho_{0,h}= \lambda/2\mu$, $r_h$ and $C_h$ are given in \eqref{r+-mu1=mu2} and \eqref{C-h-mu1=mu2}, respectively, and $\pi_0$ is given in \eqref{eqn:pi0-mu1=mu2} with $\mu_{12}=0$.
Therefore, the expected number of customers in the system is given by
\begin{equation}\label{eqn:E[L]-mu1=mu2}
  E[L_{0,h}] = \lim_{\mu_{12} \to 0} E[L_{D,h}] =  \frac{2 \rho_{0,h}}{1-\rho_{0,h}^2},
\end{equation}
where $E[L_{D,h}]$  is given in \eqref{ELD}, and the expected number of customers waiting in the queue is given by
\begin{equation}\label{eqn:E[Q]-mu1=mu2}
  E[Q_{0,h}] = \lim_{\mu_{12} \to 0} E[Q_{D,h}] =   \frac{2 \rho_{0,h}^3}{1-\rho_{0,h}^2},
\end{equation}
where $E[Q_{D,h}]$  is given in \eqref{EQD}.
\end{corollary}

\proof All above results can be obtained through elementary calculations.
\pend

\begin{remark}
The limiting results given in Corollary~\ref{cor:to-MM2-mu1=mu2} are consistent with the corresponding results for the standard $M/M/2$ model (see, for example, equations (i), (ii) and (iii) in \cite{Singh:1970}).
\end{remark}

%
%
%

\subsubsection{When $\mu_1 \neq \mu_2$}

In this case, two servers are heterogeneous. In modelling, we need to split state 1 into two states: (1,0) and (0,1) for both $M/M_D/2$ model (with dependent service times) and $M/M/2$ model (with two independent heterogeneous servers). For convenience and also without loss of generality, we assume, in both models, that
\[
    \mu_1 > \mu_2,
\]
and when both servers are idle, the next customer will join server 1 to receive its service. In the $M/M_D/2$ model, the above assumption corresponds to $p=1$ (and then $q=0$).

With the above assumption, the stationary probabilities $P_k$ for the independent model are given by (for example, see equations (7)--(10) in \cite{Singh:1970}):
\begin{align} \label{eqn:Pi-mu-1>mu_2-0}
  P_0 = & \frac{1}{1+C^*}, \\
  P_{(1,0)} = & \frac{1+\rho_0}{1+2\rho_0} \frac{\lambda}{\mu_1} P_0,\\
  P_{(0,1)} = & \frac{\rho_0}{1+2\rho_0} \frac{\lambda}{\mu_2} P_0,\\
   P_1 = & P_{(1,0)} + P_{(0,1)} = \frac{\lambda(\lambda+\mu_2)}{(1+2\rho_0)\mu_1\mu_2}P_0 = (1-\rho_0) C^* P_0, \label{eqn:Pi-mu-1>mu_2-1}  \\
  P_k = & \rho_0^{k-1} P_1, \quad k \geq 2, \label{eqn:Pi-mu-1>mu_2-k}
\end{align}
where $\rho_0 = \lambda/(\mu_1+\mu_2)$, and
\begin{equation}\label{eqn:C*}
  C^* =  \frac{\lambda(\lambda+\mu_2)}{(1+2\rho_0)\mu_1\mu_2} \frac{1}{1-\rho_0}.
\end{equation}
The expected number $E[L_0]$ of customers in the system and the expected number $E[Q_0]$ of customers waiting in the queue are given, respectively, by
\begin{align}\label{eqn:E[L]}
  E[L_0] = & \frac{1}{(1-\rho_0)^2} P_1, \\
  E[Q_0]= & \frac{\rho_0^2}{(1-\rho_0)^2} P_1, \label{eqn:E[Q]}
\end{align}
where $P_1$ is given in \eqref{eqn:Pi-mu-1>mu_2-1}.

\begin{corollary} \label{cor:to-MM2}
For the $M/M_D/2$ model with heterogeneous servers (or $\mu_1 \neq \mu_2$), as $\mu_{12} \to 0$, we have the following:
\begin{align}\label{eqn:mu{12}-to-0}
 & \lim_{\mu_{12} \to 0} r = \rho_0 = \frac{\lambda}{\mu_1+\mu_2}, \\
 & P_{[\mu_{12} = 0]} = \lim_{\mu_{12} \to 0} P(r) = \mu_1 \mu_2 (1+2\rho_0), \\
 & B_0 = \lim_{\mu_{12} \to 0} B = \lambda + \mu_2, \\
 & C_0 = \lim_{\mu_{12} \to 0} C = \left [\frac{1}{\mu_1+\mu_2} \frac{P_{[\mu_{12}=0]}}{B_0} + \frac{\rho_0}{1-\rho_0} \right ]^{-1}, \\
\end{align}
where $r$,  $P(r)$, $B$ and $C$ are given in \eqref{eqn:r+}, \eqref{eqn:P(r)}, \eqref{eqn:B} and\eqref{eqn:C}, respectively. Based on the above expressions, the stationary probabilities can be expressed in terms of $\pi_0$ as:
\begin{align}\label{eqn:pik-in-terms-of-pi0}
  \pi_0 = & C_0 \frac{\rho_0}{\lambda} \frac{P_{[\mu_{12} = 0]}}{B_0},\\
  \pi_1 = & C_0 \rho_0 = \frac{\lambda(\lambda+ \mu_2)}{\mu_1\mu_2(1+2\rho_0)} \pi_0, \\
  \pi_k = & \pi_1 \rho_0^{k-1}, \quad k \geq 2.
\end{align}
\end{corollary}

\proof All above results can be obtained through elementary calculations.
\pend

\begin{remark}
It follows from Corollary~\ref{cor:to-MM2} and elementary calculations that the expressions of $\pi_k$ for $k \geq 0$ are consistent with $P_k$ for $k \geq 0$ given in \eqref{eqn:Pi-mu-1>mu_2-0}, \eqref{eqn:Pi-mu-1>mu_2-1} and \eqref{eqn:Pi-mu-1>mu_2-k}. Therefore, $\lim_{\mu_{12} \to 0} E[L_D]$ and $\lim_{\mu_{12} \to 0}E[Q_D]$, where $E[L_D]$ and $E[Q_D]$ are given, respectively, in \eqref{ELD-1} and \eqref{EQD-1}, are consistent with the expressions given in \eqref{eqn:E[L]} and \eqref{eqn:E[Q]}, respectively. 
\end{remark}

\subsection{Solution when $\mu_1+\mu_2 \to 0$}

This is another extreme case, in which only simultaneous departure occurs. This model corresponds to the $M/M/1$ queue with bulk service of size 2, denoted by $M/M^2/1$. The infinitesimal generator of this model is given by
\[
        Q_b = \bordermatrix[{[]}]{ & 0 & 1 &2 &3 & 4 & 5 & \cdots & \cr
    0 & -\lambda & \lambda & \cr
    1 & \mu_{12} & -(\lambda+\mu_{12}) & \lambda \cr
    2 & \mu_{12} & 0 & -(\lambda+\mu_{12}) & \lambda \cr
    3 & & \mu_{12} & 0 & -(\lambda+\mu_{12}) & \lambda \cr
    4 & & & \mu_{12} & 0 & -(\lambda+\mu_{12})  & \ddots  \cr
    5 & & & & \mu_{12} & 0 & \ddots  &  \cr
    6 & & & & & \mu_{12} &  \ddots &  \cr
     \vdots & & & & && \ddots &  }.
\] 

\begin{corollary}
Under the stability condiction $\rho_b = \lambda/ (2 \mu_{12}) <1$, the stationary probability vector of this Markov chain $Q_b$ is given by
\begin{align}\label{eqn:pi-b}
  \pi_k = & C_b r_b^k, \quad k \geq 0, 
\end{align}
where
\begin{align}\label{eqn:pi-b}
  r_b = & \frac{- \mu_{12} + \sqrt{\mu_{12}^2 + 4 \lambda \mu_{12}}}{2 \mu_{12}} = \frac{-1 + \sqrt{1+ 4 \lambda/\mu_{12}}}{2}, \\
  C_b = & 1-r_b. 
\end{align}
The expected number $E[L_b]$ of customers in the system and the expected number $E[Q_b]$ of customers waiting in the queue are given, respectively, by
\begin{align}\label{eqn:E[L-b]}
  E[L_b] = & C_b \cdot \frac{r_b}{(1-r_b)^2},  \\
  E[Q_b]= & C_b \cdot \frac{r_b^3}{(1-r_b)^2}. \label{eqn:E[Q-b]}
\end{align}

\end{corollary}

\subsection{Comparisons of between $M/M_D/2$, $M/M/2$, and $M/M^2/1$}

In this section, we provide numerical comparisons of three models --- the $M/M_D/2$, $M/M/2$, and $M/M^2/1$ queues --- focusing on system performance. In particular, we examine the impact of service time dependence on the expected number of customers waiting in the system. For this purpose, we consider the following parameter values:

\begin{description}
    \item[(1)] For making the comparison meaningful, we use the same traffic intensity (system load, or server utilization) in all three models, or
    \[
        \rho = \rho_0 = \rho_b,
    \]
    where $\rho$, $\rho_0$ and $\rho_b$ are the traffic intensity for the $M/M_D/2$, $M/M/2$ and $M/M^2/1$ queues, respectively;
    
    \item[(2)] Choose three values of the traffic intensity $\rho$, representing a light, medium, and high system load, and make the arrival rate the same (for convenience of comparisons) for the three models:       
        \[
            \begin{array}{c|ccc}
            \rho & 0.1 & 0.5 & 0.9 \\ 
            \lambda & 1.5 & 7.5 & 13.5 \\ \hline
            \lambda/\rho& 15 & 15 & 15
            \end{array}
        \]
        
    \item[(3)] We calculate $E[Q_b]$, $E[Q_D]$, $E[Q_{D,h}]$, $E[Q_0]$, and $E[Q_{0,h}]$ according to equations \eqref{eqn:E[Q-b]}, 
        \eqref{EQD-1}, 
        \eqref{EQD}, 
        \eqref{eqn:E[Q]}, 
        and \eqref{eqn:E[Q]-mu1=mu2}, 
        respectively. Notice that $\rho$ is the ratio of the mean drift to the right to the mean drift to the left of the Markov chain. According to the above choice, the mean drift to the left is $\lambda/\rho =15$. The following table provides details for different cases:
       \[
            \begin{array}{ccl}
            \text{In the formula for} & \text{eqn \#} & \text{drift to left} \\
            E[Q_b]  & \eqref{eqn:E[Q-b]} & 2 \mu_{12} = 15 \\ 
            E[Q_D]  & \eqref{EQD-1} & \mu_1+\mu_2+ 2\mu_{12}=15 \\ 
            E[Q_{D,h}]  & \eqref{EQD} & 2 (\mu+\mu_{12}) =15\\ 
            E[Q_0]  & \eqref{eqn:E[Q]} & \mu_1+\mu_2=15 \\ 
            E[Q_{0,h}]  & \eqref{eqn:E[Q]-mu1=mu2} & 2 \mu =15 \\ 
            \end{array}
        \]
        With the choice of the above parameter values, for each pair of $(\rho, \lambda)$ values, $E[Q_b]$, $E[Q_0]$ and $E[Q_{0,h}]$ are completely determined. 
        
    \item[(4)] For the case of dependent service times, or $\mu_1=\mu_2=\mu$, since $\mu$ can be expressed by $\mu_{12}$ as $\mu = 7.5 - \mu_{12}$, $E[Q_{D,h}]$ is now a function of $\mu_{12}$. By choosing values of $\mu_{12}$ in its range $[0,7.5]$, impact of dependence can be analyzed. 
        
        Now, $E[Q_{D,h}]$ is completely determined, as a function of $\mu_{12}$.
        
    \item[(5)] For the case of of dependent service times with $\mu_1 \neq \mu_2$, to write $E[Q_D]$ as a function of $\mu_{12}$, we need further to specify the ratio $\theta$ of $\mu_1$ to $\mu_2$, or $\theta = \mu_1/\mu_2$, such that for fixed ratio $\theta$,
        \[
            \mu_2 = \frac{\mu_1+\mu_2}{\theta + 1}=\frac{15-2\mu_{12}}{\theta + 1},\quad \mu_1 = \frac{\theta (15-2\mu_{12})}{\theta +1}.
        \]
        Now, for each choice of $\theta$ values, we can numerically see the trend of $E[Q_D]$ as $\mu_{12}$ changes.
\end{description}

\subsection{Numerical analysis when $\mu_1 = \mu_2$}

\begin{figure}[htbp]
    \centering
      \includegraphics[width=0.7\textwidth]{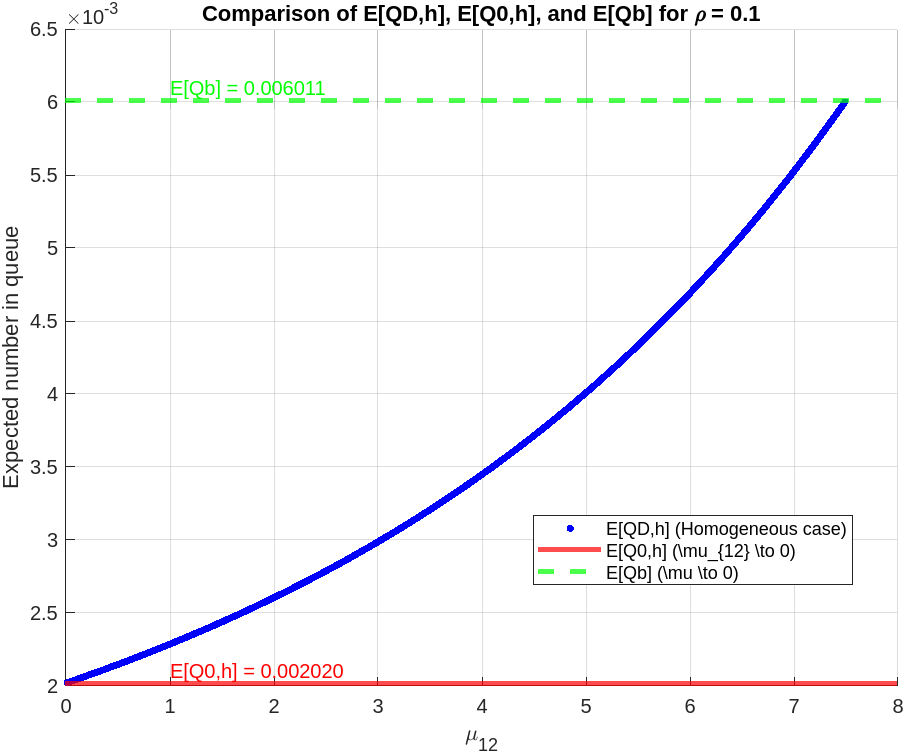}
    \caption{$\rho =0.1$ for impact on the system performance by $\mu_{12}$ when $\mu_1=\mu_2=\mu$}
     \label{fig:1-0.1}
\end{figure}
\begin{figure}[htbp]
    \centering
      \includegraphics[width=0.7\textwidth]{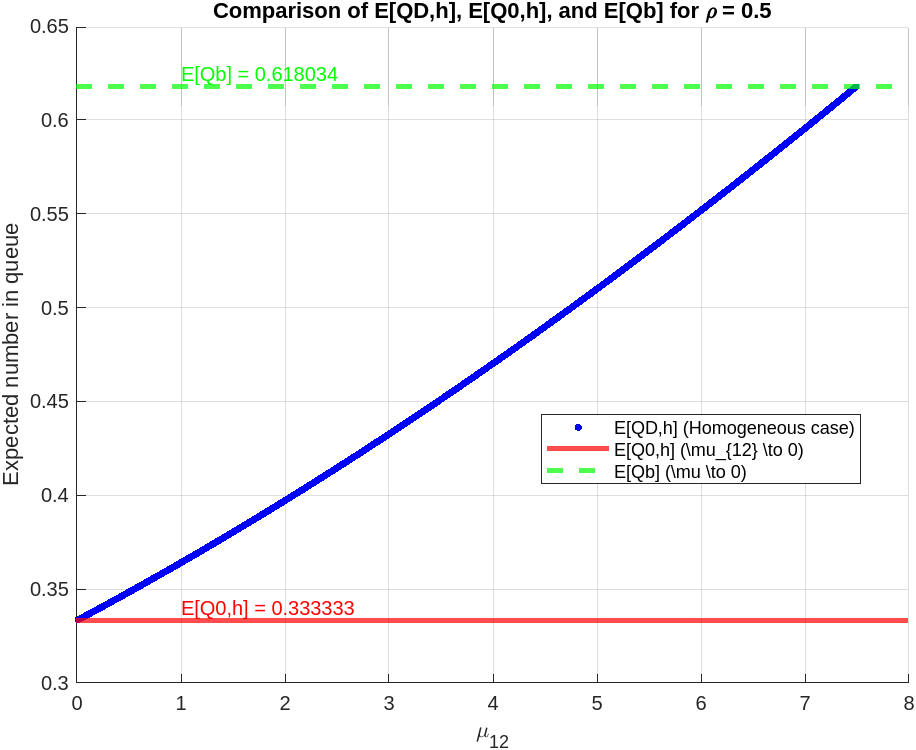}
    \caption{$\rho =0.5$ for impact on the system performance by $\mu_{12}$ when $\mu_1=\mu_2=\mu$}
     \label{fig:1-0.5}
\end{figure}
\begin{figure}[htbp]
    \centering
      \includegraphics[width=0.7\textwidth]{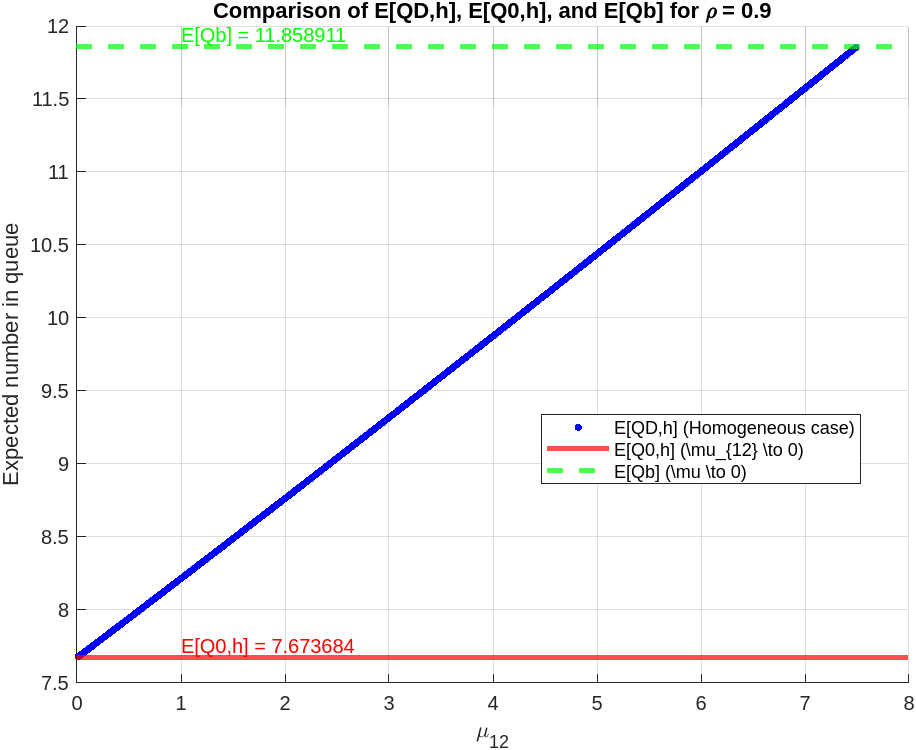}
    \caption{$\rho =0.9$ for impact on the system performance by $\mu_{12}$ when $\mu_1=\mu_2=\mu$}
     \label{fig:1-0.9}
\end{figure}

When $\mu_1=\mu_2=\mu$, the expected number $E[Q_{D,h}]$, given in \eqref{EQD}, of customers waiting in the queue is 
plotted in Figures~\ref{fig:1-0.1} -- \ref{fig:1-0.9} for $\rho = 0.1$, 0.5 and 0.9, respectively.
We also plot the two limiting cases $E[Q_{0,h}]$ and $E[Q_b]$, as horizontal lines corresponding to \( \mu_{12} \to 0 \) and \( \mu \to 0 \), respectively.
To ensure a fair comparison, the system load (or traffic intensity) is kept the same across all three models. Specifically, for each value of \( \rho \), we set \( \lambda = 15\rho \) in all three models.
It is then easy to see that the range for \( \mu_{12} \) is \([0, 7.5]\), where \( \mu_{12} = 0 \) and \( \mu_{12} = 7.5 \) correspond to the two limiting cases, respectively.

%

It is evident that the dependence between service times has a significant impact on system performance. This impact becomes increasingly pronounced as $\mu$ goes to 0. For a small value of $\rho$ (say $\rho =0.1$), the increase in $E[Q_{D,h}]$ is initially slower than linear, 
specifically, it grows more slowly than the slope of the straight line connecting \( (0, E[Q_{0,h}]) \) and \( (7.5, E[Q_b]) \) as \( \mu_{12} \) increases to a certain point. Beyond that point, the growth becomes faster than linear. 
As \( \rho \) increases, the curve of \( E[Q_{D,h}] \) approaches that of the straight line (see, for example, the curve for \( \rho = 0.9 \)).
For \( \rho = 0.1 \), 0.5, and 0.9, the maximum value of \( E[Q_{D,h}] \) is approximately 2.9757, 1.8541, and 1.5454 times as large as \( E[Q_{0,h}] \), respectively, for the model with independent servers.
These results show that the relative impact of service time dependence decreases as \( \rho \) increases—from 2.9757 to 1.5454 as \( \rho \) increases from 0.1 to 0.9.

We conjecture that, in the homogeneous case, \( E[Q_{D,h}] \) is an increasing function of \( \mu_{12} \) for any fixed value of \( \rho \).

\subsection{Numerical analysis when $\mu_1 \neq \mu_2$}

\begin{figure}[htbp]
    \centering
       \includegraphics[width=0.7\textwidth]{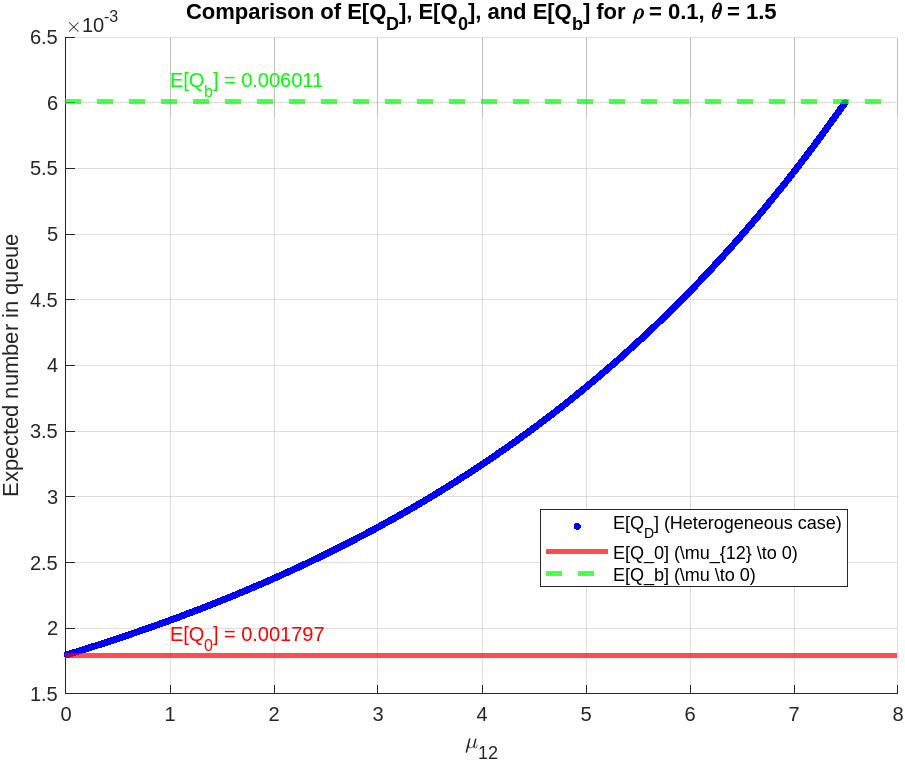}
    \caption{$\rho =0.1$ for impact on the system performance by $\mu_{12}$ when $\mu_1 \neq \mu_2$ with $\theta =1.5$}    
     \label{fig:2-0.1}
\end{figure}
\begin{figure}[htbp]
    \centering
       \includegraphics[width=0.7\textwidth]{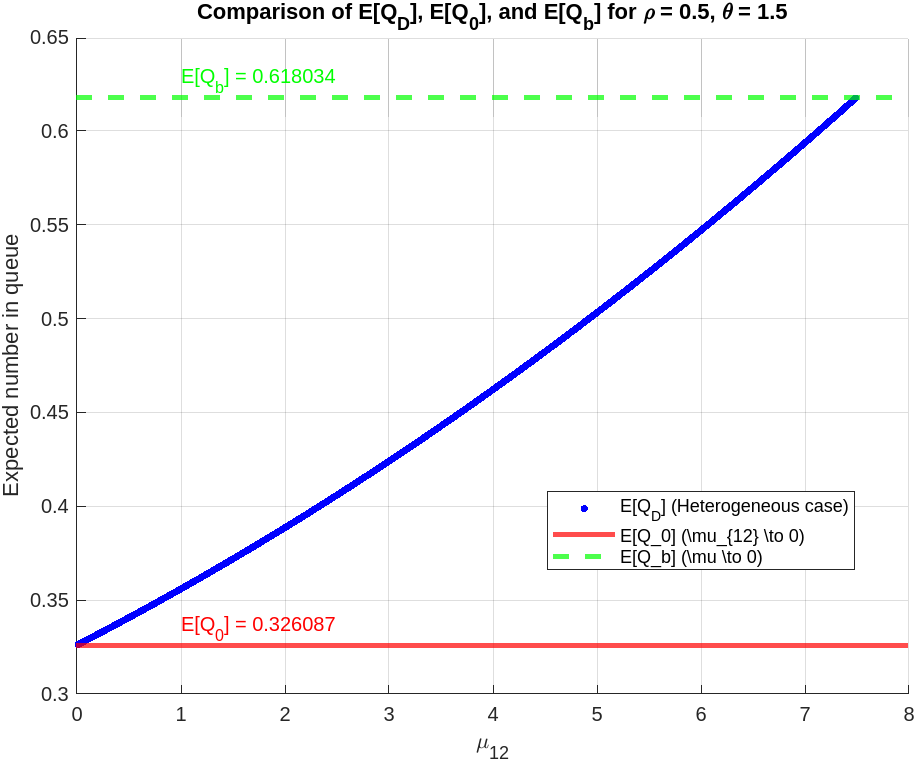}
    \caption{$\rho =0.5$ for impact on the system performance by $\mu_{12}$ when $\mu_1 \neq \mu_2$ with $\theta =1.5$}    
     \label{fig:2-0.5}
\end{figure}
\begin{figure}[htbp]
    \centering
       \includegraphics[width=0.7\textwidth]{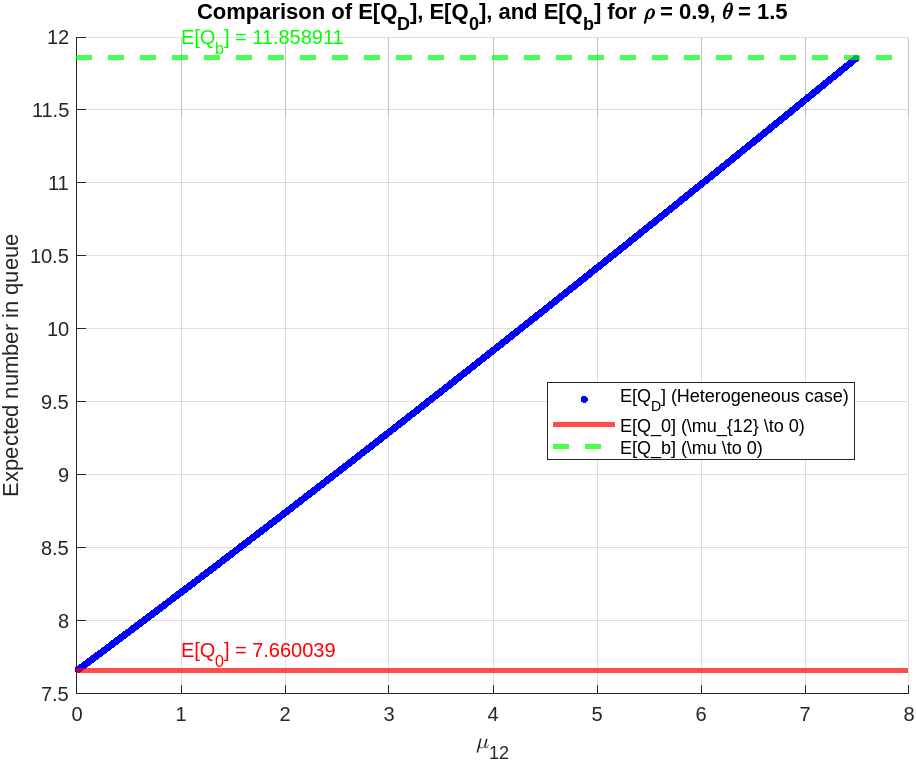}
    \caption{$\rho =0.9$ for impact on the system performance by $\mu_{12}$ when $\mu_1 \neq \mu_2$ with $\theta =1.5$}    
     \label{fig:2-0.9}
\end{figure}
\begin{figure}[htbp]
    \centering
       \includegraphics[width=0.7\textwidth]{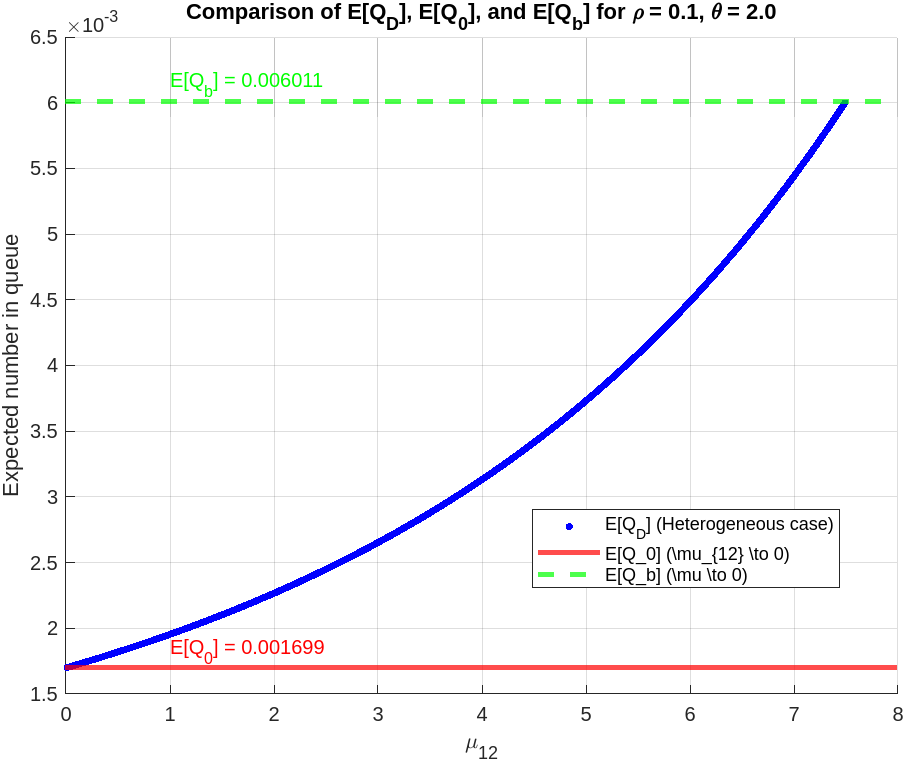}
    \caption{$\rho =0.1$ for impact on the system performance by $\mu_{12}$ when $\mu_1 \neq \mu_2$ with $\theta =2$}    
     \label{fig:3-0.1}
\end{figure}
\begin{figure}[htbp]
    \centering
       \includegraphics[width=0.7\textwidth]{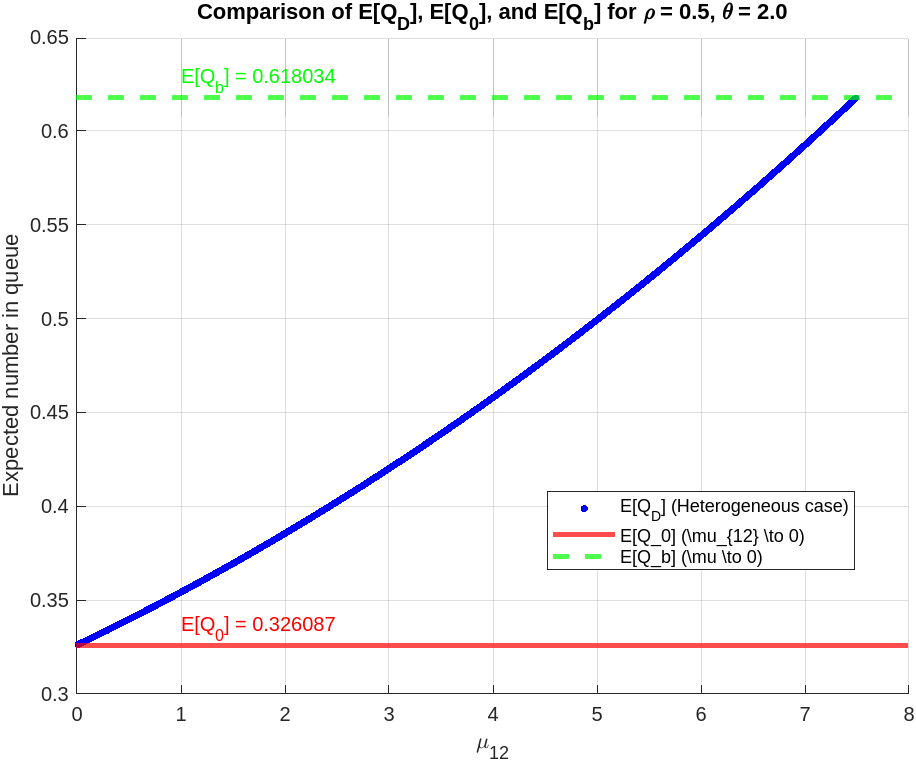}
    \caption{$\rho =0.5$ for impact on the system performance by $\mu_{12}$ when $\mu_1 \neq \mu_2$ with $\theta =2$}    
     \label{fig:3-0.5}
\end{figure}
\begin{figure}[htbp]
    \centering
       \includegraphics[width=0.7\textwidth]{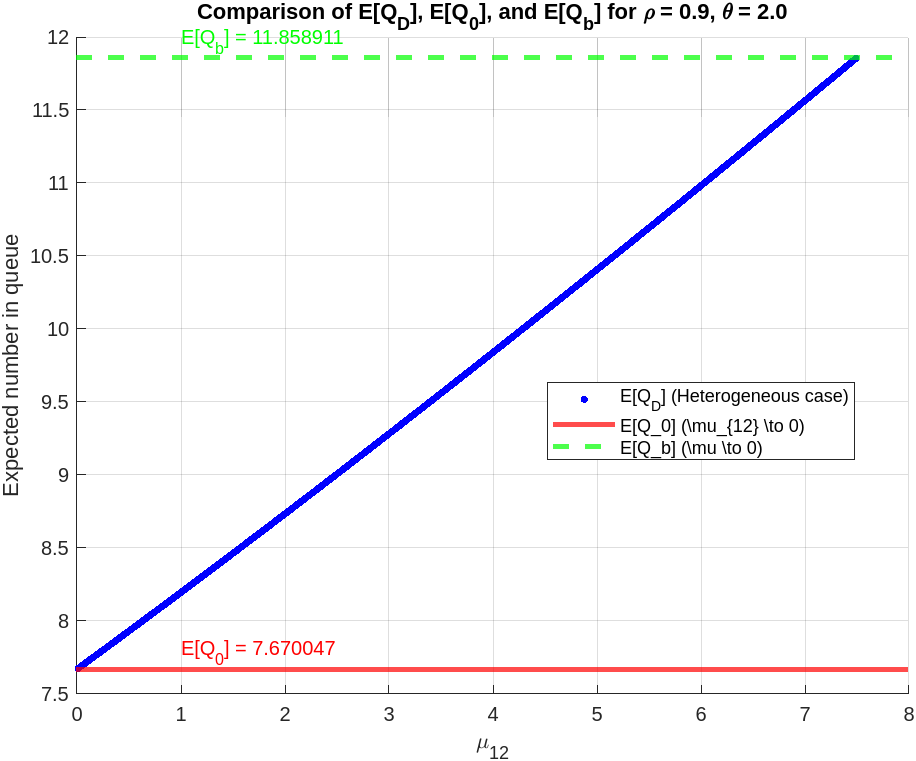}
    \caption{$\rho =0.9$ for impact on the system performance by $\mu_{12}$ when $\mu_1 \neq \mu_2$ with $\theta =2$}    
     \label{fig:3-0.9}
\end{figure}
\begin{figure}[htbp]
    \centering
       \includegraphics[width=0.7\textwidth]{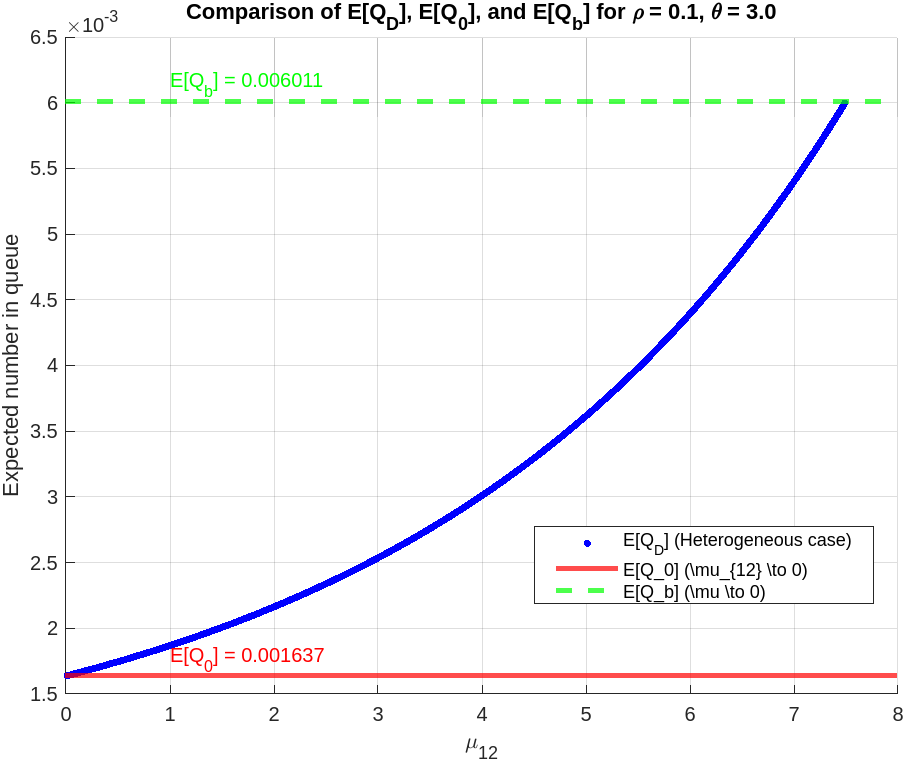}
    \caption{$\rho =0.1$ for impact on the system performance by $\mu_{12}$ when $\mu_1 \neq \mu_2$ with $\theta =3$}    
     \label{fig:4-0.1}
\end{figure}
\begin{figure}[htbp]
    \centering
       \includegraphics[width=0.7\textwidth]{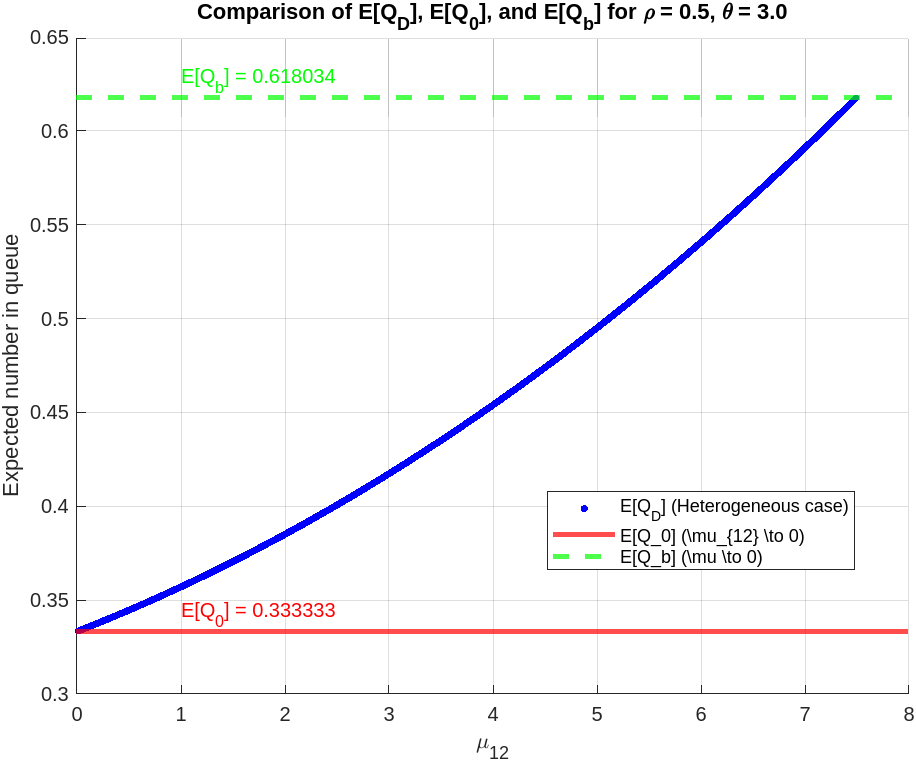}
    \caption{$\rho =0.5$ for impact on the system performance by $\mu_{12}$ when $\mu_1 \neq \mu_2$ with $\theta =3$}    
     \label{fig:4-0.5}
\end{figure}
\begin{figure}[htbp]
    \centering
       \includegraphics[width=0.7\textwidth]{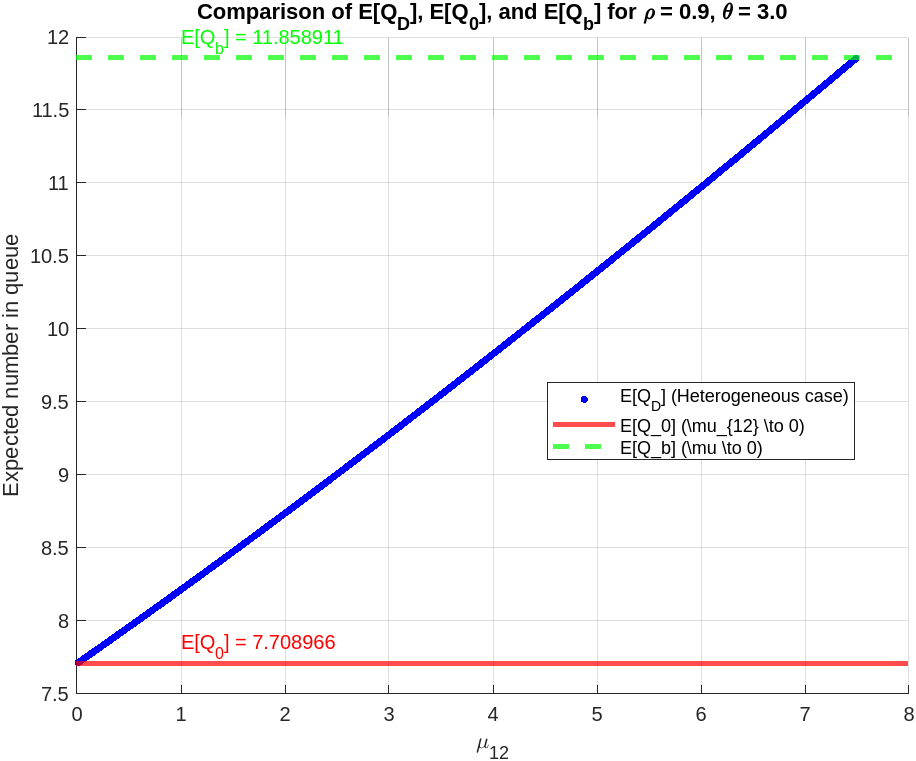}
    \caption{$\rho =0.9$ for impact on the system performance by $\mu_{12}$ when $\mu_1 \neq \mu_2$ with $\theta =3$}    
     \label{fig:4-0.9}
\end{figure}

For the heterogeneous case, where \( \mu_1 \neq \mu_2 \), we must specify the ratio \( \theta = \mu_1 / \mu_2 \), or equivalently \( \mu_1 = \theta \mu_2 \), in order to produce a specific plot of \( E[Q_D] \).
For \( \theta = 1.5 \), 2, and 3, 
Figures~\ref{fig:2-0.1} -- \ref{fig:4-0.9}
show that the impact of correlated service times is similar to that in the homogeneous server case, as long as the ratio \( \theta \) is not too large (e.g., for the three values of \( \theta \) specified above).

\begin{figure}[htbp]
    \centering
       \includegraphics[width=0.7\textwidth]{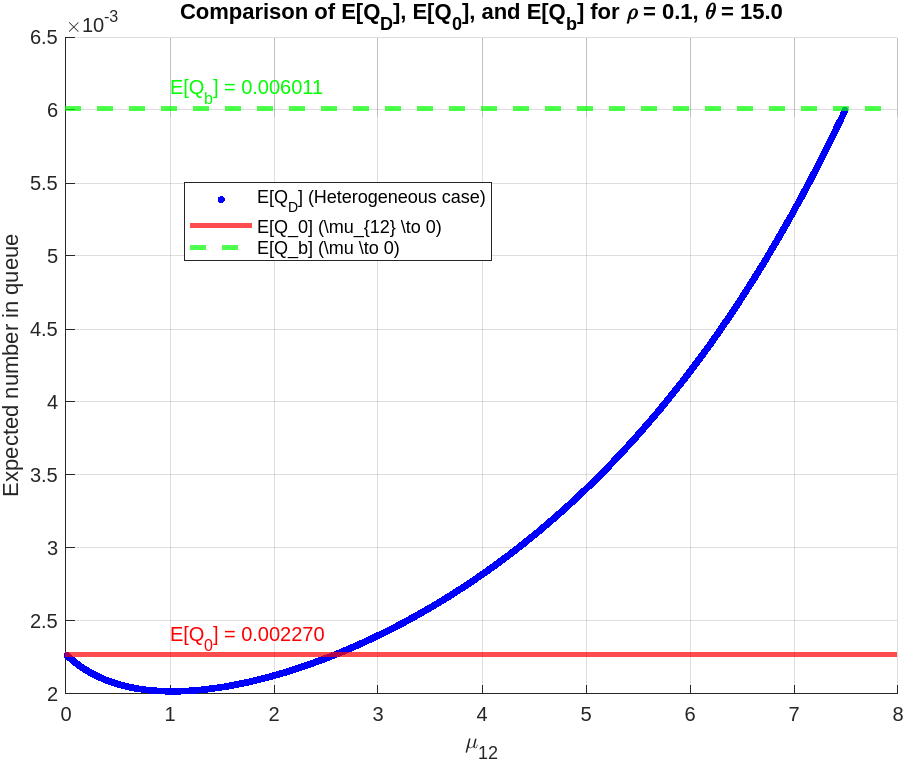}
    \caption{$\rho =0.1$ for impact on the system performance by $\mu_{12}$ when $\mu_1 \neq \mu_2$ with $\theta =15$}    
     \label{fig:5-0.1}
\end{figure}

\begin{figure}[htbp]
    \centering
       \includegraphics[width=0.7\textwidth]{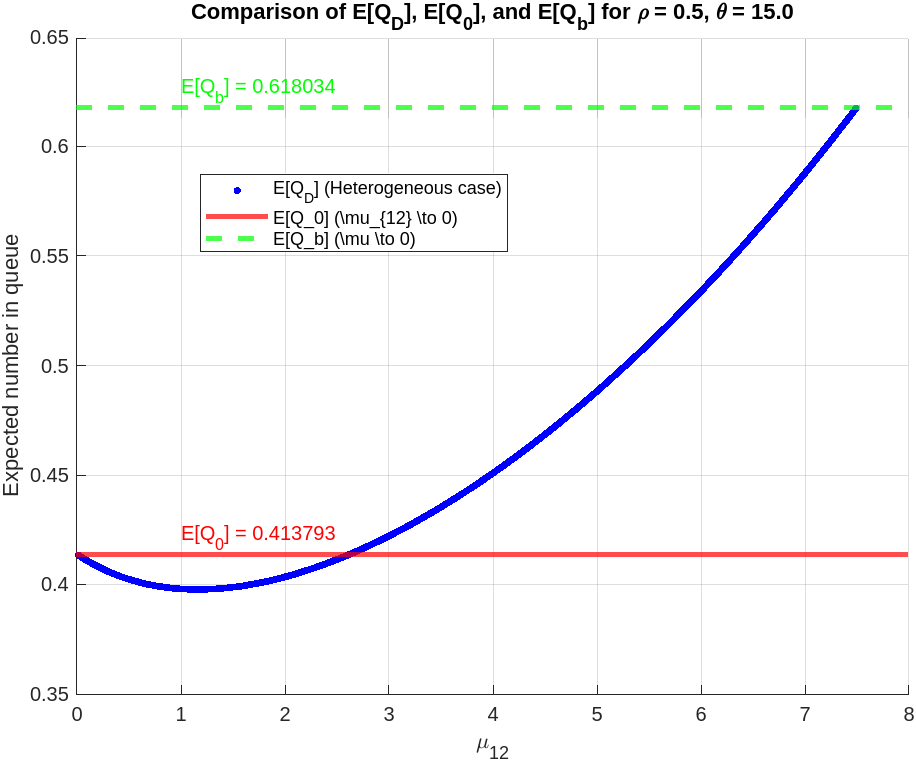}
    \caption{$\rho =0.5$ for impact on the system performance by $\mu_{12}$ when $\mu_1 \neq \mu_2$ with $\theta =15$}    
     \label{fig:5-0.5}
\end{figure}

However, It is interesting to observe that, when the dependence parameter $\theta$
is large (for example, $\theta =15$ in Figures~\ref{fig:5-0.1} and \ref{fig:5-0.5}) while the joint service rate $\mu_{12}$ remains relatively small, the expected queue length in the dependent system can be smaller than that
in the independent system, which is not intuitively obvious.

This phenomenon can be attributed to several interacting effects. First,
in this regime the dependence structure primarily induces synchronization
between the two servers rather than frequent simultaneous service
completions. As a result, the system spends less time in asymmetric
states such as $(1,0)$ or $(0,1)$, which are inefficient due to one idle
server, thereby improving effective utilization.

Second, even infrequent joint service events introduce occasional
simultaneous departures, which accelerate the return of the system to
low-congestion states (in particular, the empty state). This leads to
shorter busy periods and contributes to a reduction in the overall queue
length.

Finally, since the stationary distribution exhibits a geometric tail,
a slight reduction in the effective tail parameter (induced by the
dependence structure) results in a noticeable compression of the tail,
which further reduces $E[Q_D]$.

When $\mu_{12}$ increases, however, the effect of simultaneous service
events becomes dominant, and the performance of the dependent system
gradually deviates from that of the independent case.

\begin{figure}[htbp]
    \centering
       \includegraphics[width=0.7\textwidth]{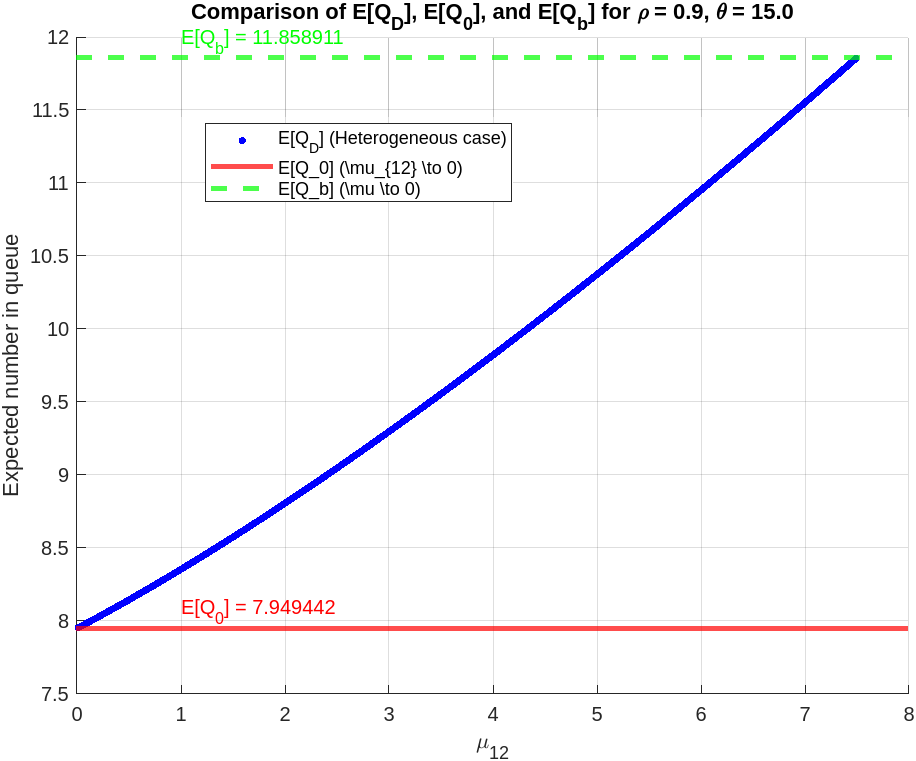}
    \caption{$\rho =0.9$ for impact on the system performance by $\mu_{12}$ when $\mu_1 \neq \mu_2$ with $\theta =15$}    
     \label{fig:5-0.9}
\end{figure}

This phenomenon becomes less noticeable as \( \rho \) increases. When \( \rho \) is sufficiently large, we no longer expect better performance from dependence; that is, \( E[Q_D] \) remains greater than \( E[Q_0] \) for all \( \mu_{12} \). When \( \rho \) is very large (e.g., \( \rho = 0.9 \)), the curve of \( E[Q_D] \) becomes nearly linear with respect to \( \mu_{12} \), even when \( \theta \) is large (as shown in Figure~\ref{fig:5-0.9}).

%

\section{Concluding remarks} \label{sec:6}

 In this paper, we proposed a new approach that leads to a Markovian queueing model for systems with positively correlated servers. This approach successfully addresses the concerns raised by Mitchell \textit{et al.}~\cite{Mitchell et al:1977}, who wrote:
\begin{quote}
“It is not clear that the birth–death equation approach can be modified to incorporate dependent service times. Moreover, any such formulation would very likely be analytically intractable.”
\end{quote}

We have shown that the birth–death equation approach can indeed be extended to handle dependent service times. However, the resulting Markov chain for the queueing process is no longer a birth–death process due to singularities in the service time density. Nevertheless, for the \( M/M_D/2 \) model, the Markov chain remains analytically tractable.

We expect that this approach can be applied to other two-dimensional queueing systems with dependent servers, and can be also extended to systems of dimension higher than two.


\begin{thebibliography}{99}

\bibitem{Avi-Itzhak-Levy:2001}
Avi-Itzhak, B. and Levy, H. (2001)
Bufferr equirement and server ordering in a tandem queue with correlated service times, 
\textit{Mathematics of Operations Research}, 
\textbf{26(2)}, 358--74.


\bibitem{Boxma:1979}
Boxma, O.J. (1979)
On a tandem queueing model with identical service times at both counters, I
\textit{Advances in applied probability},
\textbf{11(3)}, 616--643.

\bibitem{Boxma:1979b}
Boxma, O.J. (1979)
On a tandem queueing model with identical service times at both counters, II
\textit{Advances in applied probability},
\textbf{11(3)}, 644--659.

\bibitem{Boxma-Deng:2000}
Boxma, O.J. and Deng, Q. (2000)
Asymptotic behavior of the tandem queueing system with
identical service times at both queues,
\textit{Mathematical Methods of Operations Research},
\textbf{52}, 307--23.

\bibitem{Browning:1998}
Browning, S.G. (1998)
Tandem queues with blocking: A comparison between dependent and independent service,
\textit{Operations Research}, \textbf{46(3)}, 424--429.

\bibitem{Calo:1979}
Calo, S.B. (1979)
The message channel, a tandem interconnection of queues, Part I, 
\textit{IBM Report RC 6868}.

\bibitem{Calo:1980}
Calo, S.B. (1980)
The message channel, a tandem interconnection of queues, Part II, 
\textit{IBM Report RC 7170}.



\bibitem{Erlang:1909}
Erlang, A.K. (1909)
The theory of probabilities and telephone conversations,
\textit{Ny Tidsskrift for Matematik, B}, 
\textbf{20}, 33--39.

\bibitem{Choo-Conolly:1980}
Hoon Choo, Q. and  Conolly, B. (1980)
Waiting time analysis for a tandem queue with correlated service,
\textit{European journal of operational research},
\textbf{4(5)}, 337--345.

\bibitem{Gelenbe:1989}
Gelenbe, E. (1989)
Random neural networks with negative and positive signals and product form solution, 
\textit{Neural Computation}, 
\textbf{1},  502--510.

\bibitem{Gelenbe:1991}
Gelenbe, E. (1991)
Product-form queueing networks with negative and positive customers, 
\textit{Journal of Applied Probability}, \textbf{28}, 656--663.

\bibitem{Gumbel:1960}
Gumbel, H. (1960)
Waiting lines with heterogeneous servers,
\textit{Operations Research}, 
\textbf{8}, 504--511.

%

\bibitem{Kelly:1976}
Kelly, F.P. (1976)
Networks of queues,
\textit{Adv. Appl. Prob.}, \textbf{8}, 416--432.

\bibitem{Kelly:1979}
Kelly, F.P. (1979)
\textit{Reversibility and Stochastic Networks},
Wiley, Chichester.

\bibitem{Kelly:1982}
Kelly, F.P (1982)
The throughput of a series of buffers,
\textit{Adv. in Appl. Probab.}, \textbf{14}, 633--653.

\bibitem{Kleinrock:1964}
Kleinrock, L. (1964)
\textit{Communication Nets: Stochastic Message Flow and Delay},
McGraw-Hill, New York.

\bibitem{Kotz-et at:2000}
Kotz, S.,  Balakrishnan, N., and Johnson, N.L. (2000)
\textit{Continuous Multivariate Distributions: Models and Applications},
John Wiley \& Sons.

\bibitem{Lin-Dou-Kuriki:2019}
Lin, G.D., Dou, X., and Kuriki, S. (2019)
The bivariate lack-of-memory distributions,
\textit{Sankhyā: The Indian Journal of Statistics}, 
\textbf{81}, 273--297. 

\bibitem{Marshall-Olkin:1967a}
Marshall, A.W. and Olkin, I. (1967)
A multivariate exponential distribution,
\textit{Journal of the American Statistical Association}.
\textbf{62}, 30--44.



\bibitem{Mitchell et al:1977}
Mitchell, C.R., Paulson, A.S. and Beswick, C.A. (1977)
The effect of correlated exponential service times on single server tandem queues,
\textit{Naval Research Logistics Quarterly},
\textbf{24(1)}, 95--112.



\bibitem{Pang-Whitt:2012}
Pang, G. and Whitt, W. (2012)
Infinite-server queues with batch arrivals and dependent service times,
\textit{Probability in the Engineering and Informational Sciences},
\textbf{26}, 197--220.

\bibitem{Pang-Whitt:2012-MSOM}
Pang, G. and Whitt, W. (2012)
The impact of dependent service times on large-scale service systems,
\textit{Manufacturing and Service Operations Management},
\textbf{14(2)}, 226--278.

\bibitem{Pang-Whitt:2013}
Pang, G. and Whitt, W. (2013)
Two-parameter heavy-traffic limits for infinite-server queues with dependent service times, \textit{Queueing Systems},
\textbf{73(2)}, 119--146.


\bibitem{Pinedo-Wolff:1982}
Pinedo, M. and Wolff, R.W. (1982)
A comparison between tandem queues with dependent and independent service times,
\textit{Operations Research}, \textbf{30(3)}, 464--479.

\bibitem{Sandmann:2007}
Sandmann, W. (2007)
Performance evaluation of dependent two-stage services, 
\textit{Proceedings of the 21st European Conference on Modelling and Simulation (ECMS)},
74--9.

\bibitem{Sandmann:2010}
Sandmann, W. (2010)
Delays in a series of queues: Independent versus identical service times,
\textit{The IEEE symposium on Computers and Communications}, 32--37.


\bibitem{Sandmann:2012}
Sandmann, W. (2012)
Delays in a series of queues with correlated service times,
\textit{Journal of Network and Computer Applications},
\textbf{35(5)}, 1415--1423.

\bibitem{Singh:1970}
Singh, V.P. (1970)
Two-server Markovian queues with balking: heterogeneous vs. homogeneous servers,
\textit{Operations Research}, 
\textbf{18}, 145--159.

\bibitem{Wolff:1982}
Wolff, R.W. (1982)
Tandem queues with dependant service times in light traffic,
\textit{Operations research}, \textbf{30}, 619--635.

\bibitem{Ziedins:1993}
Ziedins, I. (1993)
Tandem queues with correlate service times and finite capacity,
\textit{Math. O.R.}, \textbf{18}, 901--915.

\end{thebibliography}
\end{document}